\documentclass[bj,authoryear]{imsart}
\usepackage{eurosym}
\usepackage{amssymb,amsmath,amsfonts,amsthm,enumerate}
\usepackage{graphicx,subcaption,float,relsize}
\usepackage{natbib,url,rotating}
\usepackage{hyperref}
\usepackage{sectsty}
\usepackage[justification=centering]{caption}

\startlocaldefs
\numberwithin{equation}{section}
\theoremstyle{plain}

\allowdisplaybreaks
\newtheorem {theorem}{Theorem}[section]
\newtheorem {assumption}{Assumption}

\newtheorem{example}[theorem]{Example}
\newtheorem{lemma}[theorem]{Lemma}
\newtheorem{remark}{Remark}

\endlocaldefs
\begin{document}

\begin{frontmatter}
\title{A bootstrapped test of covariance stationarity based on orthonormal transformations}
\runtitle{A bootstrapped test of covariance stationarity}

\begin{aug}
\author[A]{Jonathan B. Hill}
\author[B]{Tianqi Li}
\address[A]{Dept. of Economics, University of North Carolina}
\address[B]{Dept. of Economics, University of North Carolina}

\end{aug}

\begin{abstract}
We propose a covariance stationarity test for an otherwise
dependent and possibly globally non-stationary time series. We work in a generalized version of the new setting in \cite{JinWangWang2015}, who  exploit \cite{Walsh1923} functions in order to compare sub-sample covariances with the
full sample counterpart. They impose strict stationarity under the null,
only consider linear processes under either hypothesis in order to achieve a parametric estimator for an inverted high dimensional asymptotic covariance matrix, and do not consider any other orthonormal basis. 
Conversely, we work with a general orthonormal basis under mild conditions that include Haar wavelet and Walsh functions,
and we allow for linear or nonlinear processes with possibly non-iid innovations. This is important in macroeconomics and finance where nonlinear feedback and random volatility occur in many settings. We completely sidestep asymptotic covariance matrix estimation and inversion by bootstrapping a max-correlation difference statistic, where the maximum is taken over the correlation lag $h$ and basis generated sub-sample counter $k$ (the number of systematic
samples). We achieve a higher feasible rate of increase for the maximum lag
and counter $\mathcal{H}_{T}$ and $\mathcal{K}_{T}$. Of particular note, our test is capable of detecting breaks in variance, and distant, or very mild, deviations from stationarity.
\end{abstract}


\end{frontmatter}


\section{Introduction\label{sec:intro}}

Assume $\{X_{t}$ $:$ $t$ $\in $ $\mathbb{Z}\}$ is a possibly non-stationary
time series process in $\mathcal{L}_{2}$. We want to test whether $X_{t}$ is
covariance stationary, without explicitly assuming stationarity under the
null hypothesis, allowing for linear or nonlinear processes with a possibly
non-iid innovation, and a general memory property. Such generality is
important in macroeconomics and finance where nonlinear feedback and non-iid
innovations occur in many settings due to asymmetries and random volatility,
including exchange rates, bonds, interest rates, commodities, and asset
return levels and volatility. Popular models for such time series include
symmetric and asymmetric GARCH, Stochastic Volatility, nonlinear ARMA-GARCH,
and switching models like smooth transition autoregression. See, e.g., \cite%
{Terasvirta1994}, \cite{Gray1996} and \cite{FrancqZakoian2019}.

Evidence for nonstationarity, whether generally or in the variance or
autocovariances, has been suggested for many economic time series, where
breaks in variance and model parameters are well known %
\citep[e.g.][]{BusettTaylor2003,Perron2006,HendryMassmann2007,GianettoRaissi2015}%
. Knowing whether a time series is globally nonstationary has large
implications for how analysts approach estimation and inference. Indeed, it
effects whether conventional parametric and semi-(non)parametric model
specifications are correct. Pretesting for deviations from global
stationarity therefore has important practical value.

There are many tests in the literature on covariance stationarity, and
concerning locally stationary processes. Tests for stationarity based on
spectral or second order dependence properties have a long history, where
pioneering work is due to \cite{PriestleySubbaRao1969}. Spectrum-based tests
with $\mathcal{L}_{2}$-distance components have many versions. \cite%
{Paparoditis_2010} uses a rolling window method to compare subsample local
periodograms against a full sample version. The maximum is taken over the $%
\mathcal{L}_{2}$-distance between periodograms over all time points. An
asymptotic theory for the max-statistic, however, is not provided, although
an approximation theory is (see their Lemmas 1 and 3). Furthermore,
conforming with many offerings in the literature, under the null $X_{t}$ is
a linear process with iid Gaussian innovations. \cite{Dette_etal_2011} study
locally stationary processes, and impose linearity with iid Gaussian
innovations. Their statistic is based on the minimum $\mathcal{L}_{2}$%
-distance between a spectral density and its version under stationarity, and
local power is non-trivial against $T^{1/4}$-alternatives. \cite%
{Aue_etal2009} propose a nonparametric test for a break in covariance for
multivariate time series based on a version of a cumulative sum statistic.

Wavelet methods have arisen in various forms recently. \cite%
{vonSachsNeumann2000}, using technical wavelet decomposition components from 
\cite{NeumannVonSachs1997}, propose a Haar wavelet based localized
periodogram test of covariance stationarity for locally stationary processes %
\citep[cf.][]{Dahlhaus1997,Dahlhaus2009}, but neglect to characterize power.
Haar wavelet functions form an orthonormal basis on $\mathcal{L}_{2}[0,1)$,
but the proposed frequency domain tests are complicated, a local power
analysis is not feasible, and empirical power may be weak (see simulation
evidence from \cite{JinWangWang2015}).

\cite{DwivediSubbaRao2011} and \cite{JentschSubbaRao2015} use the discrete
Fourier transform [DFT] $J_{T}(\omega _{k})$ $=$ $(2\pi
T)^{-1/2}\sum_{t=1}^{T}X_{t}\exp \left\{ it\omega _{k}\right\} $ at
canonical frequencies $\omega _{k}$ $=$ $2\pi k/T$ and $1$ $\leq $ $k$ $\leq 
$ $T$. \cite{DwivediSubbaRao2011} generate a portmanteau statistic from a
normalized sample DFT covariance, exploiting the fact that an uncorrelated
DFT implies second order stationarity. \cite{Nason2013} presents a
covariance stationarity test based on Haar wavelet coefficients of the
wavelet periodogram, they assume linear local stationarity, and do not treat
local power. See also \cite{NasonVonSachsKroisandt2000}.

In a promising offering in the wavelet literature, \cite{JinWangWang2015}
[JWW] exploit so-called \textit{Walsh functions} (akin to \textquotedblleft
global square waves\textquotedblright\ although not truly wavelets; cf. \cite%
{Walsh1923} and their implied systematic samples for comparing sub-sample
covariances with the full sample one. They utilize a sample-size dependent
maximum lag $\mathcal{H}_{T}$ and maximum systematic sample counter $%
\mathcal{K}_{T}$, and show their Wald test exhibits non-negligible local
power against $\sqrt{T}$-alternatives. They do not consider any other
orthonormal transformation because Walsh functions, they argue, have
\textquotedblleft \textit{desirable properties}\textquotedblright\ based
primarily on simulation evidence, asymptotic independence of a sub-sample
and sample covariance difference $(\sqrt{T}(\hat{\gamma}_{h}^{(k_{1})}$ $-$ $%
\hat{\gamma}_{h}),\sqrt{T}(\hat{\gamma}_{h}^{(k_{2})}$ $-$ $\hat{\gamma}%
_{h}))$ across systematic samples $k_{1}$ $\neq $ $k_{2}$, and joint
asymptotic normality (JWW, p. 897). It seems, however, that such theoretical
properties are available irrespective of the orthonormal basis used,
although we do not provide a proof. See Section \ref{sec: Walsh}, below, for
definitions and notation. We do, however, find in the sequel that the Walsh
basis has superlative properties vis-\`{a}-vis a Haar wavelet basis.

JWW's asymptotic analysis is driven by local stationarity and linearity $%
X_{t}=\sum_{i=0}^{\infty }\psi _{i}Z_{t-i}$, with zero mean iid $Z_{t}$, and 
$E|Z_{t}|^{4+\delta }$ $<$ $\infty $, $\delta $ $>$ $0$, which expedites
characterizing a parametric asymptotic covariance matrix estimator. The iid
and linearity assumptions, however, rule out many important processes,
including nonlinear models like regime switching, and random coefficient
processes, and any process with a non-iid error (e.g. nonlinear ARMA-GARCH).
JWW's Wald-type test statistic requires an inverted parametric variance
estimator that itself requires five tuning parameters and choice of two
kernels.\footnote{%
One tuning parameter $\lambda $ $\in $ $(0,.5)$ governs the number $Q_{T}$ $%
= $ $[T^{\lambda }]$ of sample covariances that enter the asymptote
covariance matrix estimator (see their p. 899). The remaining four $%
(c_{1},c_{2};\xi_{1},\xi _{2})$\ are used for kernel bandwidths $b_{j}$ $=$ $%
c_{j}T^{-\xi_{j}}$, $j$ $=$ $1,2$, for computing the kurtosis of the iid
process $Z_{t}$ under linearity (see p. 902-903). The authors set $c_{j}$
equal to $1.2$ times a so-called "\textit{crude scale estimate}" which is
nowhere defined.} Indeed, most of the tuning parameters only make sense
under linearity given how they approach asymptotic covariance matrix
estimation.

Now define the lag $h$ autocovariance coefficient at time $t$:%
\begin{equation*}
\gamma _{h}(t)\equiv E\left[ \left( X_{t}-E\left[ X_{t}\right] \right)
\left( X_{t-h}-E\left[ X_{t-h}\right] \right) \right] \text{, }h=0,1,...
\end{equation*}%
The hypotheses are:%
\begin{eqnarray}
&&H_{0}:\gamma _{h}(s)=\gamma _{h}(t)=\gamma _{h}\text{ }\forall s,t\text{, }%
\forall h=0,1,...\text{ (cov. stationary)}  \label{null} \\
&&H_{1}:\gamma _{h}(s)\neq \gamma _{h}(t)\text{ for some }s\neq t\text{ and }%
h=0,1,...\text{ (cov. nonstationary).}  \notag
\end{eqnarray}%
Under $H_{0}$ $X_{t}$ is second order stationary, and the alternative is 
\textit{any} deviation from the null: the autocovariance differs across time
at some lag, allowing for a (lag zero) break in variance. The null
hypothesis otherwise accepts the possibility of global nonstationarity.

In this paper we do away with parametric assumptions on $X_{t}$, and impose
either a mixing or physical dependence property that allows us to bound the
number of usable covariance lags $\mathcal{H}_{T}$ and systematic samples $%
\mathcal{K}_{T}$. The conditions allow for global nonstationarity under
either hypothesis, allowing us to focus the null hypothesis only on second
order stationarity. We show that use of the physical dependence construct in 
\cite{Wu2005} is a boon for bounding $\mathcal{H}_{T}$ since it allows for
slower than geometric memory decay, and ultimately need only hold uniformly
over ($h,k$). The mixing condition imposed, however, requires geometric
decay, and must ultimately hold \textit{jointly} over all lags $h$ $=$ $%
1,...,\mathcal{H}_{T}$ ($h$ is the covariance lag, and $k$ is a systematic
subsample counter). The latter leads to a much diminished upper bound on the
maximum lag growth $\mathcal{H}_{T}$ $\rightarrow $ $\infty $. This may be
of independent interest given the recent rise of high dimensional central
limit theorems under weak dependence %
\citep[eg.][]{ZhangWu2017,ChangJiangShao2023,ChangChenWu2024}.

Rather than operate on a Wald statistic constructed from a specific
orthonormal transformation of covariances, our statistic is the maximum
generic orthonormal transformed sample correlation coefficient, where the
maximum is taken over ($h,k$) with increasing upper bounds $(\mathcal{H}_{T},%
\mathcal{K}_{T})$. By working in a generic setting we are able to make
direct comparisons, and combine bases for possible power improvements.

We provide examples of Haar wavelet and Walsh functions in Sections \ref%
{sec: Walsh} and \ref{sec:Haar}, and show how they yield different
systematic samples. This suggests a power improvement may be available by
using \textit{multiple} orthonormal transforms. As JWW (p. 897) note,
however, clearly other orthonormal transformations are feasible, although
simulation evidence agrees with their suggestion that the Walsh basis works
quite well.

We use a dependent wild bootstrap for the resulting test statistic, allowing
us to sidestep asymptotic covariance matrix estimation, a challenge
considering we do not assume a parametric form, and the null hypothesis
requires us to look over a large set of $(h,k)$. We sidestep all of JWW's
tuning parameters, and require just one governing the block size for the
bootstrap. We ultimately achieve a significantly better upper bound on the
rate of increase for $(\mathcal{H}_{T},\mathcal{K}_{T})$ than JWW. Penalized
and weighted versions of our test statistic are also possible, as in JWW and 
\cite{HillMotegi2020} respectively. There is, though, no compelling theory
to justify penalties on $(h,k)$ in our setting, and overall a non-penalized
and unweighted test statistic works best in practice.

Note that \cite{HillMotegi2020} study the max-correlation statistic for a
white noise test, and only show their limit theory applies for some
increasing maximum lag $\mathcal{H}_{T}$, but do not derive an upper bound.
In the present paper we use a different asymptotic theory, derive upper
bounds for $\mathcal{H}_{T}$ and $\mathcal{K}_{T}$, and of course do not
require a white noise property under $H_{0}$.

\citet[Section
2.6]{JinWangWang2015} rule out the use of autocorrelations because, they
claim, if the sample variance were included, i.e. $h$ $\geq $ $0$, then
consistency may not hold because the limit theory neglects the joint
distribution of $\hat{\gamma}_{0}$ and the correlation differences. We show
for our proposed test that the difference between full sample and systematic
sample autocorrelations at lag zero asymptotically reveals whether $%
E[X_{t}^{2}]$ is time dependent. Further, our test is consistent whether
non-stationarity is caused by variances, or covariances, or both. See
Section \ref{sec:max_corr_H1} and Example \ref{ex:var}. Our proposed test is
consistent against a general (nonparametric) alternative, and exhibits
nontrivial power against a sequence of $\sqrt{T}$-local alternatives.

The max-correlation difference is particularly adept at revealing subtle
deviations from covariance stationarity, similar to results revealed in \cite%
{HillMotegi2020}. Consider a distant form of a model treated in 
\citet[Model
I]{Paparoditis_2010} and 
\citet[Section 3.2: models NVI,
NVII]{JinWangWang2015}, $X_{t}$ $=$ $.08\cos \{1.5$ $-$ $\cos (4\pi
t/T)\}\epsilon _{t-d}$ $+$ $\epsilon _{t}$ with large $d$ (JWW use $d$ $=$ $%
1 $ or $6$). JWW's test exhibits trivial power when $d$ $\geq $ $20$, while
the max-correlation difference is able to detect this deviation from the
null even when $d$ $\geq $ $50$. The reason is the same as that provided in 
\cite{HillMotegi2020}: the max-correlation difference operates on the single
most useful statistic, while Wald and portmanteau statistics congregate many
standardized covariances that generally provide little relevant information
under a weak signal.

In Section \ref{sec:max_corr_Walsh} we develop the test statistic. Sections %
\ref{sec:max_corr_asym} and \ref{sec:bootstrap} present asymptotic theory
and the bootstrap method and theory. We then perform a Monte Carlo study in
Section \ref{sec:monte_carlo}, and conclude with Section \ref{sec:conclusion}%
. The supplemental material \cite{supp_mat_covstat_2024} contains all
proofs, an empirical study concerning international interest rates, and
complete simulation results.

We use the following notation. $[z]$ rounds $z$ to the nearest integer. $%
\mathcal{L}_{2}$ is the space of square integrable random variables; $%
\mathcal{L}_{2}[a,b)$ is the class of square integrable functions on $[a,b)$%
. $||\cdot ||_{p}$ and $||\cdot ||$ are the $L_{p}$ and $l_{2}$ norms
respectively, $p$ $\geq $ $1$. Let $\mathbb{Z}$ $\equiv $ $%
\{...-2,-1,0,1,2,...\}$, and $\mathbb{N}$ $\equiv $ $\{0,1,2,...\}$. $K$ $>$ 
$0$ is a finite constant whose value may be different in different places. $%
awp1$ denotes \textquotedblleft \textit{asymptotically with probability
approaching one}\textquotedblright . Write $\max_{\mathcal{H}_{T}}$ $=$ $%
\max_{0\leq h\leq \mathcal{H}_{T}}$. $\max_{\mathcal{K}_{T}}$ $=$ $%
\max_{1\leq k\leq \mathcal{K}_{T}}$ and $\max_{\mathcal{H}_{T},\mathcal{K}%
_{T}}$ $=$ $\max_{0\leq h\leq \mathcal{H}_{T},1\leq k\leq \mathcal{K}_{T}}$.
Similarly, $\max_{\mathcal{H}_{T}}a(h,\tilde{h})$ $=$ $\max_{0\leq h,\tilde{h%
}\leq \mathcal{H}_{T}}a(h,\tilde{h})$, etc. $|a|_{+}$ $\equiv $ $a$ $\vee $ $%
0$.

\section{Max-correlation with orthonormal transformation\label%
{sec:max_corr_Walsh}}

Our test statistic is the maximum of an orthonormal transformed sample
covariance. In order to build intuition, we first derive the test statistic
under Walsh function and Haar wavelet-based bases. We then set up a general
environment, and present the main results.

In order to reduce notation, assume here $\mu $ $\equiv $ $E[X_{t}]$ $=$ $0$
is known. In practice this is enforced by using $X_{t}$ $-$ $\bar{X}$ where $%
\bar{X}$ $\equiv $ $1/T\sum_{t=1}^{T}X_{t}$. In 
\citet[Lemmas B.3 and
B.3$^{*}$]{supp_mat_covstat_2024} we prove using $X_{t}$ $-$ $\bar{X}$ or $%
X_{t}$ $-$ $\mu $ leads to identical results asymptotically. Thus in proofs
of the main results we simply assume $\mu $ $=$ $0$.

\subsection{Walsh functions\label{sec: Walsh}}

The following class of Walsh functions $\{W_{i}(x)\}$ $\equiv $ $\{W_{i}(x)$ 
$:$ $i$ $=$ $0,1,2,...\}$ define a complete orthonormal basis in $\mathcal{L}%
_{2}[0,1)$. The functions $W_{i}(x)$ are defined recursively 
\citep[see,
e.g.,][]{Walsh1923,AhmedRao1975,Stoffer1987,Stoffer1991}:%
\begin{equation*}
W_{0}(x)=1\text{ for }x\in \lbrack 0,1)\text{ and }W_{1}(x)=\left\{ 
\begin{array}{cc}
1, & x\in \lbrack 0,.5) \\ 
-1, & x\in \lbrack .5,1)%
\end{array}%
\right. ,
\end{equation*}%
and for any $i$ $=$ $1,2,...,$%
\begin{equation*}
W_{2i}(x)=\left\{ 
\begin{array}{ll}
W_{i}(2x), & x\in \lbrack 0,.5) \\ 
\left( -1\right) ^{i}W_{i}(2x-1), & x\in \lbrack .5,1)%
\end{array}%
\right. \text{and }W_{2i+1}(x)=\left\{ 
\begin{array}{ll}
W_{i}(2x), & x\in \lbrack 0,.5) \\ 
\left( -1\right) ^{i+1}W_{i}(2x-1), & x\in \lbrack .5,1)%
\end{array}%
\right. .
\end{equation*}%
In the $\{-1,1\}$-valued\ sequence $\{W_{i}(x)$ $:$ $i$ $=$ $0,1,2,...\}$, $%
i $ indexes the number of zero crossings, yielding a square shaped
wave-form. See Figure \ref{fig:bases}, below, and see \citet[Figure
5]{Stoffer1991} and \citet[Figure
1]{JinWangWang2015} and their references. The $k^{th}$ discrete Walsh
functions used in this paper are then for $t$ $=$ $1,...,T$:%
\begin{equation*}
\left\{ \mathcal{W}_{k}(1),...,\mathcal{W}_{k}(T)\right\} \text{ where }%
\mathcal{W}_{k}(t)=W_{k}((t-1)/T).
\end{equation*}

Now define the covariance coefficient for a covariance stationary time
series, $\gamma _{h}$ $\equiv $ $E[X_{t}X_{t-h}]$, and denote the usual
(co)variance estimator $\hat{\gamma}_{h}$ $\equiv $ $1/T%
\sum_{t=1}^{T-h}X_{t}X_{t+h},$ $h$ $\in $ $\mathbb{N}.$ JWW use $\{\mathcal{W%
}_{i}(x)\}$ to construct a set of discrete Walsh covariance transformations:
for some integer $\mathcal{K}$ $\geq $ $1$,%
\begin{equation}
\hat{\gamma}_{h}^{W(k)}\equiv \frac{1}{T}\sum_{t=1}^{T-h}X_{t}X_{t+h}\left\{
1+\left( -1\right) ^{k-1}\mathcal{W}_{k}(t)\right\} ,\text{ }h=0,1,...,T-1,%
\text{ and }k=1,2,...,\mathcal{K}.  \label{gk}
\end{equation}%
As they point out, a sequence of systematic (sub)samples $\boldsymbol{T}%
_{k}^{W}$ $:$ $k$ $=$ $1,2,...,\mathcal{K}$ in the time domain can be
defined on the basis of Walsh functions:%
\begin{equation*}
\boldsymbol{T}_{k}^{W}\equiv \left\{ t\in T:\left( -1\right) ^{k-1}\mathcal{W%
}_{k}(t)=1\right\} .
\end{equation*}%
Now let $\mathcal{N}_{k}$ be the smallest power of $2$ that is at least $k.$
The first systematic sample is the first half of the sample time domain $%
\boldsymbol{T}_{1}^{W}$ $=$ $\{1,...,[T/2]\}$; the second is the middle half 
$\boldsymbol{T}_{2}^{W}$ $=$ $\{[T/4],[T/4]+1,...,[3T/4]\}$; the third $%
\boldsymbol{T}_{3}$ is the first and third time blocks, and so on. Notice $%
\boldsymbol{T}_{k}^{W}$ consists of $(k$ $+$ $1)/2$ blocks with at least $[T/%
\mathcal{N}_{k}]$ elements. Thus, when $h$ $<$ $T/\mathcal{N}_{k}$ then $%
\hat{\gamma}_{h}^{W(k)}$\ is just an estimate of $\gamma _{h}$ on the $%
k^{th} $ systematic sample:%
\begin{equation*}
\hat{\gamma}_{h}^{W(k)}=\frac{1}{T}\sum_{t=1}^{T-h}X_{t}X_{t+h}\left\{
1+\left( -1\right) ^{k-1}\mathcal{W}_{k}(t)\right\} =\frac{2}{T}\sum_{t\in 
\boldsymbol{T}_{k}}X_{t}X_{t+h}.
\end{equation*}%
The condition $h$ $<$ $T/\mathcal{N}_{k}$ holds asymptotically in the
Section \ref{sec:bootstrap} bootstrap setting.

The difference between the $k^{th}$ systematic sample and full sample
estimators is: 
\begin{equation*}
\hat{\gamma}_{h}^{W(k)}-\hat{\gamma}_{h}=\left( -1\right) ^{k-1}\frac{1}{T}%
\sum_{t=1}^{T-h}X_{t}X_{t+h}\mathcal{W}_{k}(t).
\end{equation*}%
Notice the $\{-1,1\}$-valued nature of $\mathcal{W}_{k}(t)$ yields a
sub-sample comparison: $\hat{\gamma}_{h}^{W(k)}$ $-$ $\hat{\gamma}_{h}$ $=$
\linebreak $1/T\sum_{t\in \boldsymbol{T}_{k}}X_{t}X_{t+h}$ $-$ $%
1/T\sum_{t\neq \boldsymbol{T}_{k}}X_{t}X_{t+h}$. Our test is based on the
maximum $|\hat{\gamma}_{h}^{W(k)}$ $-$ $\hat{\gamma}_{h}|$, in which case
the multiple $\left( -1\right) ^{k-1}$ is irrelevant. We therefore drop it
everywhere. Under the null hypothesis and mild assumptions this difference
is $O_{p}(1/\sqrt{T})$ at all lags $h$ and for all systematic samples $k$.
Thus, a test statistic can be constructed from $\sqrt{T}(\hat{\gamma}%
_{h}^{W(k)}$ $-$ $\hat{\gamma}_{h})$.

\subsection{Haar wavelet functions\label{sec:Haar}}

Define the usual Haar wavelet functions $\psi _{k,m}(x)$ $\equiv $ $%
2^{k/2}\psi (2^{k}x$ $-$ $m)$ with $x$ $\in $ $\mathbb{R}$, where $0$ $\leq $
$k$ $\leq $ $\mathcal{K}_{T}$ for some integer sequence $\{\mathcal{K}_{T}\}$%
, $0$ $\leq $ $m$ $\leq $ $2^{k}$ $-$ $1$, and mother wavelet %
\citep{Haar1910}:%
\begin{equation*}
\psi (x)=\left\{ 
\begin{array}{rr}
1, & \text{ \ }x\in \lbrack 0,.5) \\ 
-1, & x\in \lbrack .5,1) \\ 
0 & \text{otherwise}%
\end{array}%
\right. .
\end{equation*}%
Haar functions $\{\psi _{k,m}(x)\}$ form a complete orthonormal basis in $%
\mathcal{L}[0,1)$. The discretized version is:%
\begin{equation*}
\Psi _{k,m}(t)\equiv \psi _{k,m}(\left( t-1\right) /T)=2^{k/2}\psi
(2^{k}\left( t-1\right) /T-m).
\end{equation*}%
Systematic samples derived from $\{\Psi _{k,m}(t)\}$ are generally too
\textquotedblleft local\textquotedblright : $1/T\sum_{t=1}^{T-h}X_{t}X_{t+h}%
\Psi _{2,m}(t)$, for example, compares just the first eighth to the second
eighth subsample ($m$ $=$ $0$); the third eighth to the fourth eighth
subsample ($m$ $=$ $1$); and so on.

In order to yield a test statistic that compares sub-sample complements,
comparable to Walsh functions, we compile $(\psi _{k,m}(x),\Psi _{k,m}(t))$
over $0$ $\leq $ $m$ $\leq $ $2^{k}$ $-$ $1$. Set $\psi _{0}(x)$ $=$ $I(0$ $%
\leq $ $x$ $\leq $ $1)$, and for $k$ $=$ $0,1,...$ 
\begin{eqnarray*}
\psi _{k+1}(x) &\equiv &\frac{1}{2^{k/2}}\sum_{m=0}^{2^{k}-1}\psi
_{k,m}(x)=\sum_{m=0}^{2^{k}-1}\psi (2^{k}x-m) \\
\Psi _{k+1}(t) &\equiv &\frac{1}{2^{k/2}}\sum_{m=0}^{2^{k}-1}\Psi
_{k,m}(t)=\sum_{m=0}^{2^{k}-1}\psi (2^{k}\left( t-1\right) /T-m).
\end{eqnarray*}%
We set $\psi _{0}(x)$ $=$ $I(0$ $\leq $ $x$ $\leq $ $1)$ in order to unify
the local alternative analysis below, similar to the Walsh basis. It can be
shown that $\psi _{k}(x)$ $\in $ $\{-1,1\}$, and $\{\psi _{k}(x)$ $:$ $1$ $%
\leq $ $k$ $\leq $ $\mathcal{K}_{T}\}$ forms a complete orthonormal basis:
see Lemma B.1 in \cite{supp_mat_covstat_2024} for this and other properties.
Now, in the same manner as (\ref{gk}), define for $k=1,2,...,\mathcal{K}$:%
\begin{equation*}
\hat{\gamma}_{h}^{H(k)}\equiv \frac{1}{T}\sum_{t=1}^{T-h}X_{t}X_{t+h}\left\{
1+\Psi _{k}(t)\right\} ,\text{ }h=0,1,...,T-1.
\end{equation*}

\begin{figure}[ht]
\caption{Orthonormal bases}
\label{fig:bases}%
\begin{subfigure}{0.35\textwidth}
        \includegraphics[width=5cm,height=7cm]{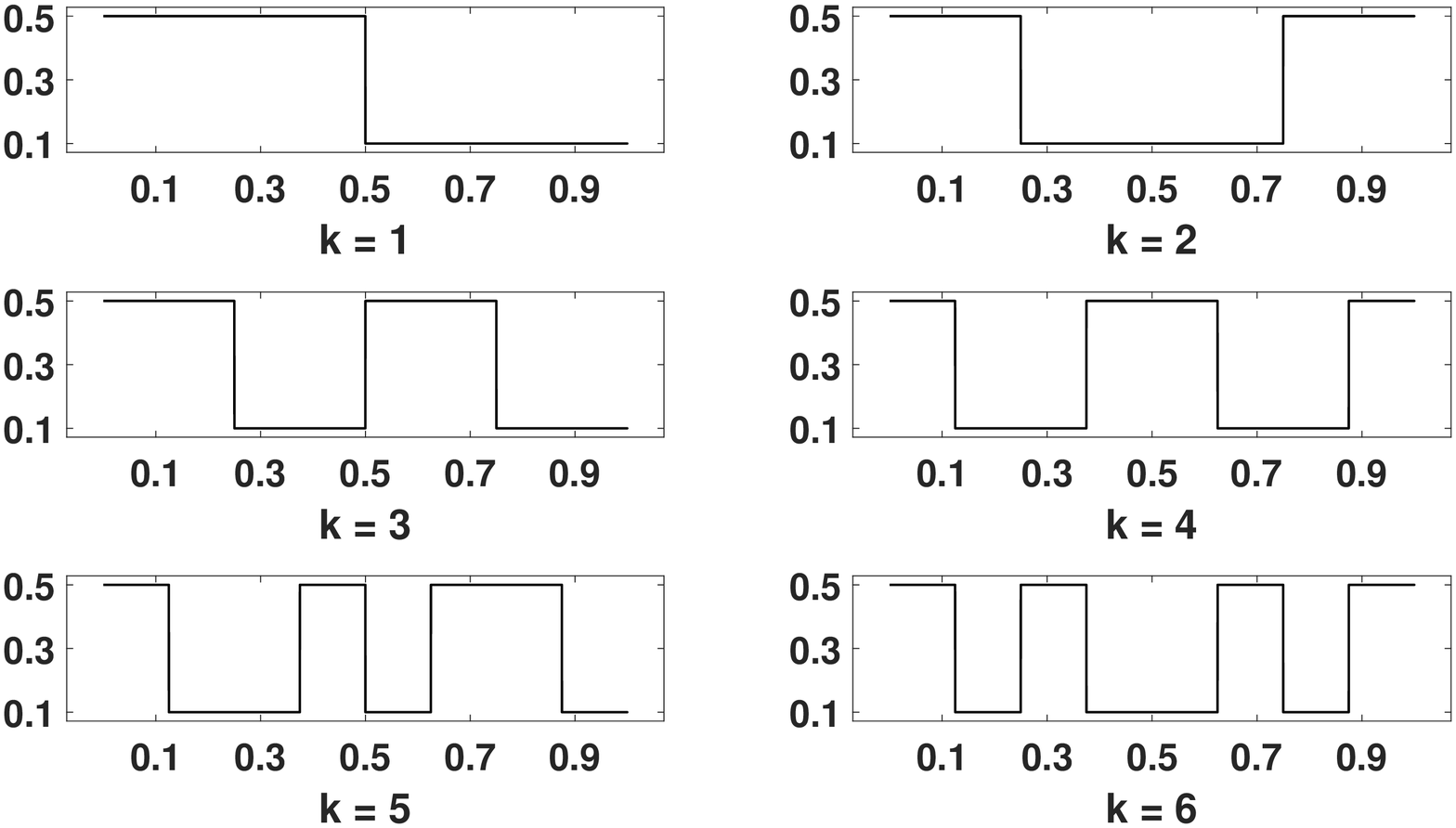}
        \caption{Walsh functions $\{W_{k}(x)\}_{k=1}^{6}$}
    \end{subfigure}\hspace{0mm} 
\begin{subfigure}{0.35\textwidth}
        \includegraphics[width=5cm,height=7cm]{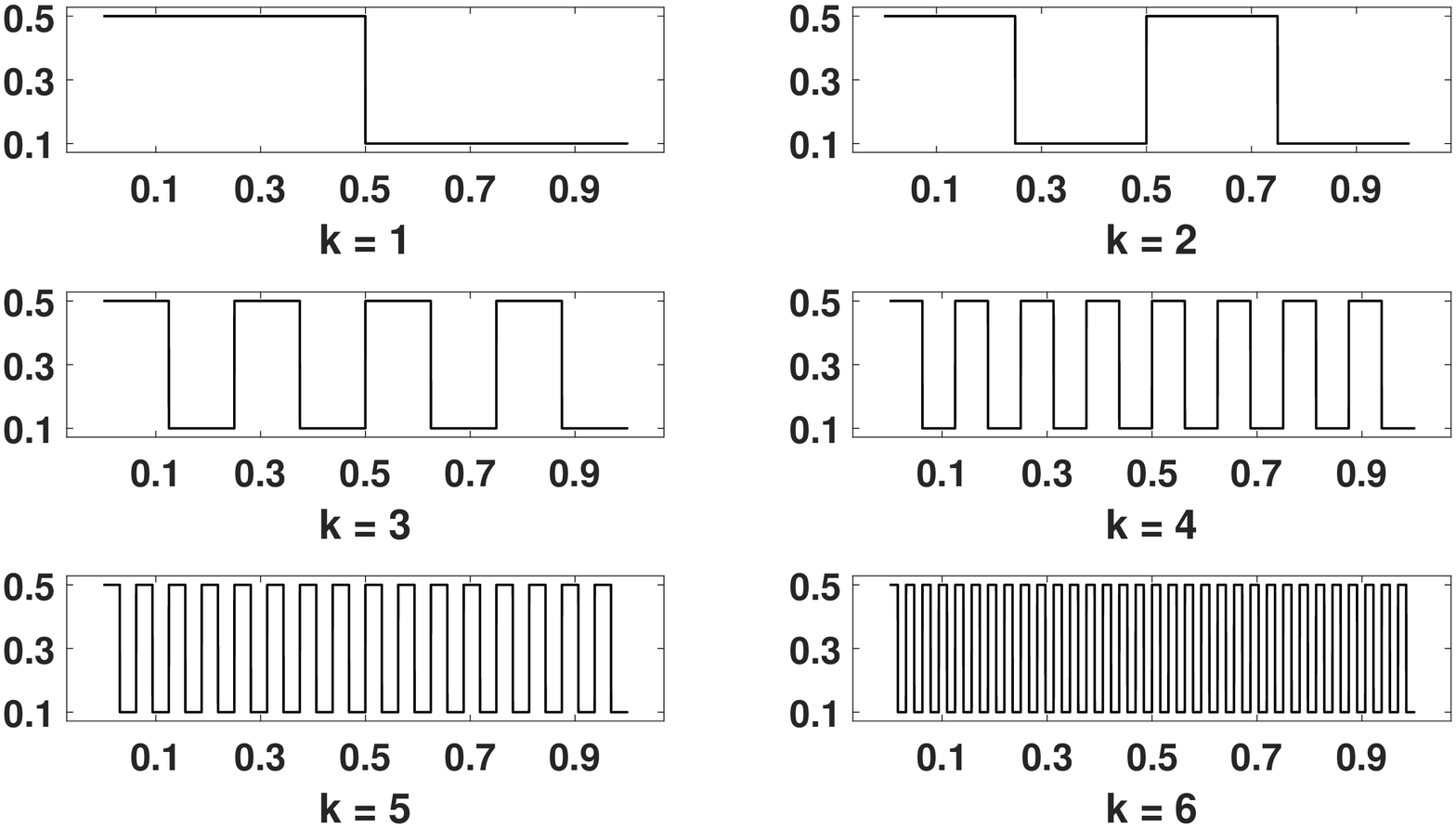}
        \caption{Composite Haar $\{\protect\psi_{k}(x)\}_{k=1}^{6} $}
    \end{subfigure}
\end{figure}

The discretized Haar functions $\Psi _{k,m}(t)$ generate systematic samples $%
\boldsymbol{T}_{k}^{H}$ $\equiv $ $\{t$ $\in $ $T$ $:$ $%
\bigcup_{m=0}^{2^{k}-1}(\Psi _{k,m}(t)$ $=$ $1)\}.$ This yields the first
sample half $\boldsymbol{T}_{0}^{H}$ $=$ $\{1,...,[T/2]\}$; the first and
third quarter subsamples $\boldsymbol{T}_{1}^{H}$ $=$ $%
\{1,...,[T/4];1+[T/2],...,[3T/4]\}$; the first, third, fifth and seventh
eights $\boldsymbol{T}_{2}^{H}$; etc. See Figure \ref{fig:bases} for plots
of Walsh and composite Haar functions $W_{k}(x)$ and $\psi _{k}(x)$, $%
k=1,...,6$.

Walsh and Haar systematic samples are quite different for $k$ $\geq $ $2$. $%
\boldsymbol{T}_{k}^{W}$ involves fewer interspersed subsample segments, in
some cases of varying lengths, while $\boldsymbol{T}_{k}^{H}$ have $2^{k}$
segments of equal length $[T/2^{k}]$ in all cases (ignoring truncation due
to the lag $h$). Haar subsamples are therefore non-redundant only when $T/2^{%
\mathcal{K}_{T}}$ $\geq $ $1$, hence we need $\mathcal{K}_{T}$ $\leq $ $\ln
(T)/\ln (2)$.

Indeed, it can be shown that the two bases coincide in the sense that $%
W_{k_{1}}(x)$ $=$ $\psi _{k_{2}}(x)$ for all $x$ and only $(k_{1},k_{2})$ $%
\in $ $\{(1,1),(3,2),(3,7)\}$, or in all other cases for $x$ on a subset of $%
[0,1)$ with Lebesgue measure $1/2$. Roughly speaking, only $50\%$ of the
data points in $\hat{\gamma}_{h}^{W(k_{1})}$ $-$ $\hat{\gamma}_{h}$ are the
same as those in $\hat{\gamma}_{h}^{H(k_{2})}$ $-$ $\hat{\gamma}_{h}$ for
nearly all systematic samples $(k_{1},k_{2})$. Thus the two bases are
intrinsically different, suggesting potential advantages and weaknesses
against certain deviations from the null.

\subsection{Max-correlation orthonormal transforms\label{maxcorrWalsh}}

Now let $\{\mathcal{B}_{k}(x)$ $:$ $0$ $\leq $ $k$ $\leq $ $\mathcal{K}\}$
denote a $\{-1,1\}$-valued orthonormal basis on $\mathcal{L}[0,1)$, $%
\mathcal{B}_{0}(x)$ $=$ $I(0$ $\leq $ $0$ $\leq $ $1),$ let $B_{k}(t)$ $%
\equiv $ $\mathcal{B}_{k}((t$ $-$ $1)/T)$ be the discretized version, and
define the usual subsample covariance for this generic discrete basis $\hat{%
\gamma}_{h}^{(k)}$ $\equiv $ $1/T\sum_{t=1}^{T-h}X_{t}X_{t+h}\{1$ $+$ $%
B_{k}(t)\}$.

Define the sample correlation coefficient:%
\begin{equation*}
\hat{\rho}_{h}\equiv \frac{\hat{\gamma}_{h}}{\hat{\gamma}_{0}},
\end{equation*}%
and a set of discrete orthonormal correlation transformations, over
systematic sample $k$: 
\begin{equation*}
\hat{\rho}_{h,1}^{(k)}\equiv \frac{\hat{\gamma}_{h}^{(k)}}{\hat{\gamma}_{0}}=%
\frac{1}{\hat{\gamma}_{0}}\times \frac{1}{T}\sum_{t=1}^{T-h}X_{t}X_{t+h}%
\left\{ 1+B_{k}(t)\right\} ,\text{ }k=1,2,...,\mathcal{K}.
\end{equation*}%
Thus, the difference between systematic sample and full sample estimators is:%
\begin{equation}
\hat{\rho}_{h,1}^{(k)}-\hat{\rho}_{h}=\frac{1}{\hat{\gamma}_{0}}\frac{1}{T}%
\sum_{t=1}^{T-h}X_{t}X_{t+h}B_{k}(t)=\frac{\hat{\gamma}_{h}^{(k)}-\hat{\gamma%
}_{h}}{\hat{\gamma}_{0}}.  \label{rho1}
\end{equation}

The correlation difference $\hat{\rho}_{h,1}^{(k)}$ $-$ $\hat{\rho}_{h}$ is
sensible even at lag $0$, considering 
\begin{equation*}
\hat{\rho}_{0,1}^{(k)}-\hat{\rho}_{0}=\frac{\hat{\gamma}_{0}^{(k)}}{\hat{%
\gamma}_{0}}-1.
\end{equation*}%
Thus, under nonstationarity $\hat{\rho}_{0,1}^{(k)}$ $\overset{p}{%
\nrightarrow }$ $1$ for some systematic sample $k$\ when $\hat{\gamma}%
_{0}^{(k)}/\hat{\gamma}_{0}$ $\overset{p}{\nrightarrow }$ $1$; that is, when
the second moment $E[X_{t}^{2}]$ is not constant over $t$.

Alternatively, we may incorporate\ the systematic sample variance estimators 
$\hat{\gamma}_{0}^{(k)}$. The autocorrelation estimator in that case
becomes, for example: 
\begin{eqnarray*}
\hat{\rho}_{h,2}^{(k)} &\equiv &\frac{1}{\hat{\gamma}_{0}}\times \frac{1}{T}%
\sum_{t=1}^{T-h}X_{t}X_{t+h}+\frac{1}{\hat{\gamma}_{0}^{(k)}}\frac{1}{T}%
\sum_{t=1}^{T-h}X_{t}X_{t+h}B_{k}(t) \\
&=&\hat{\rho}_{h}+\frac{1}{\hat{\gamma}_{0}^{(k)}}\frac{1}{T}%
\sum_{t=1}^{T-h}X_{t}X_{t+h}B_{k}(t),
\end{eqnarray*}%
hence%
\begin{equation*}
\hat{\rho}_{h,2}^{(k)}-\hat{\rho}_{h}=\frac{\frac{1}{T}%
\sum_{t=1}^{T-h}X_{t}X_{t+h}B_{k}(t)}{\frac{1}{T}\sum_{t=1}^{T}X_{t}^{2}%
\left\{ 1+B_{k}(t)\right\} }.
\end{equation*}%
At lag $0$ notice:%
\begin{equation*}
\hat{\rho}_{0,2}^{(k)}-\hat{\rho}_{0}=\frac{1}{\hat{\gamma}_{0}^{(k)}}\frac{1%
}{T}\sum_{t=1}^{T}X_{t}^{2}B_{k}(t)=\frac{\hat{\gamma}_{0}^{(k)}-\hat{\gamma}%
_{0}}{\hat{\gamma}_{0}^{(k)}}=1-\frac{\hat{\gamma}_{0}}{\hat{\gamma}%
_{0}^{(k)}}.
\end{equation*}%
Compare this to $\hat{\rho}_{0,1}^{(k)}$ $-$ $\hat{\rho}_{0}$ $=$ $\hat{%
\gamma}_{0}^{(k)}/\hat{\gamma}_{0}$ $-$ $1$. Thus, again $\hat{\rho}%
_{0,2}^{(k)}$ $\overset{p}{\nrightarrow }$ $1$ for some systematic sample $k$
when $E[X_{t}^{2}]$ is not constant over $t$.

Asymptotically $\hat{\rho}_{h,i}^{(k)}$ are identical in probability. Indeed,%
\begin{eqnarray}
\hat{\rho}_{h,1}^{(k)}-\hat{\rho}_{h,2}^{(k)} &=&\left( \frac{1}{\hat{\gamma}%
_{0}}-\frac{1}{\hat{\gamma}_{0}^{(k)}}\right) \frac{1}{T}%
\sum_{t=1}^{T-h}X_{t}X_{t+h}B_{k}(t)  \notag \\
&=&\left( \hat{\gamma}_{0}^{(k)}-\hat{\gamma}_{0}\right) \frac{1}{\hat{\gamma%
}_{0}\hat{\gamma}_{0}^{(k)}}\left( \hat{\gamma}_{h}^{(k)}-\hat{\gamma}%
_{h}\right) .  \label{rho1_rho2}
\end{eqnarray}%
This reveals $\hat{\rho}_{h,1}^{(k)}-\hat{\rho}_{h,2}^{(k)}$ for each $h$ $%
\geq $ $0$ simultaneously captures systematic sample differences in variance 
\textit{and} covariance. Under $H_{0}$ and general conditions presented in
Section \ref{sec:max_corr_asym}, $\max_{\mathcal{H}_{T},\mathcal{K}_{T}}$ $|%
\hat{\gamma}_{h}^{(k)}$ $-$ $\hat{\gamma}_{h}|$ and $|\hat{\gamma}_{0}$ $-$ $%
\gamma _{0}|$ are $O_{p}(1/\sqrt{T})$, where $\{\mathcal{H}_{T},\mathcal{K}%
_{T}\}$ are sequences defined below with $\mathcal{H}_{T}$ $\rightarrow $ $%
\infty $ and $\mathcal{K}_{T}$ $\rightarrow \infty $. Thus:%
\begin{equation*}
\max_{\mathcal{H}_{T},\mathcal{K}_{T}}\left\vert \sqrt{T}\left( \hat{\rho}%
_{h,1}^{(k)}-\hat{\rho}_{h,2}^{(k)}\right) \right\vert =O_{p}(1/\sqrt{T}).
\end{equation*}%
Under $H_{1}$, however, it holds that $\sqrt{T}\max_{\mathcal{H}_{T},%
\mathcal{K}_{T}}$ $|\hat{\rho}_{h,1}^{(k)}-\hat{\rho}_{h,2}^{(k)}|$ $\overset%
{p}{\rightarrow }$ $\infty $ \textit{if and only if} $E[X_{t}^{2}]$ and $%
E[X_{t}X_{t-h}]$ \textit{for some} $h$ $\geq $ $1$ are time dependent. This
suggests $\mathcal{D}_{T}$ $\equiv $ $\max_{\mathcal{H}_{T},\mathcal{K}_{T}}$
$|\sqrt{T}(\hat{\rho}_{h,1}^{(k)}$ $-$ $\hat{\rho}_{h,2}^{(k)})|$ could be
used as a third test statistic: $\mathcal{D}_{T}$ will reject $H_{0}$
asymptotically with power approaching one when $X_{t}$ is non-stationary in
variance \textit{and} autocovariance at some lag $h$ $\geq $ $1$.
Conversely, either $\max_{\mathcal{H}_{T},\mathcal{K}_{T}}$ $|\sqrt{T}(\hat{%
\rho}_{h,i}^{(k)}$ $-$ $\hat{\rho}_{h})|$ is consistent against $H_{1}$ in
general: power is one asymptotically if $E[X_{t}^{2}]$ \textit{and/or} some $%
E[X_{t}X_{t-h}]$ are time dependent.

In order to focus ideas we only consider $\hat{\rho}_{h,1}^{(k)}$, so put:%
\begin{equation*}
\hat{\rho}_{h}^{(k)}\equiv \hat{\rho}_{h,1}^{(k)}.
\end{equation*}%
The proposed test statistic is therefore the maximum normalized $\hat{\rho}%
_{h}^{(k)}$ $-$ $\hat{\rho}_{h}$ over $(h,k)$:%
\begin{equation*}
\mathcal{M}_{T}\equiv \sqrt{T}\max_{\mathcal{H}_{T},\mathcal{K}%
_{T}}\left\vert \hat{\rho}_{h}^{(k)}-\hat{\rho}_{h}\right\vert =\frac{1}{%
\hat{\gamma}_{0}}\max_{\mathcal{H}_{T},\mathcal{K}_{T}}\left\vert \frac{1}{%
\sqrt{T}}\sum_{t=1}^{T-h}X_{t}X_{t+h}B_{k}(t)\right\vert .
\end{equation*}%
By construction $\mathcal{M}_{T}$ uses the most informative systematic
sample correlation difference. Notice we search over \textit{all} lags $h$ $%
\in $ $\{0,...,\mathcal{H}_{T}\}$.

A penalized version is also possible:%
\begin{equation*}
\mathcal{M}_{T}^{(p)}\equiv \max_{\mathcal{H}_{T},\mathcal{K}_{T}}\left\{ 
\sqrt{T}\left\vert \hat{\rho}_{h}^{(k)}-\hat{\rho}_{h}\right\vert -\mathcal{P%
}(h,k)\right\} ,
\end{equation*}%
where $\mathcal{P}(h,k)$ is a non-random, positive, strictly monotonically
increasing function of $h$ and $k$. JWW use $\mathcal{P}(h,k)$ $=$ $p_{h}$ $+
$ $q_{h}$ with AIC-like lag penalty $p_{h}$ $=$ $2h$ in an order
selection-type Wald statistic. This is sensible considering the Wald
statistic is pointwise asymptotically chi-squared with mean $2h$ for each $k$
\citep[see also][]{InglotLedwina2006}. For $q_{k}$ they use $\sqrt{k-1}$
based on empirical power considerations, and to satisfy a required technical
condition\ \citep[eq. (3.4)]{JinWangWang2015}.

In our non-Wald setting a similar reasoning for choosing $\mathcal{P}(h,k)$ $%
=$ $p_{h}$ $+$ $q_{h}$ does not apply, nor do we have any comparable
requirements for penalizing $k$. Indeed, a compelling reason for
\textquotedblleft penalizing\textquotedblright\ $\mathcal{M}_{T}$ at all
would be to counter the loss of observations at higher lags or to control
for lag specific heterogeneity, but that historically is ameliorated with
weights, for example $\max_{\mathcal{H}_{T},\mathcal{K}_{T}}\{\sqrt{T}%
\mathfrak{W}_{T,h}^{(k)}|\hat{\rho}_{h}^{(k)}$ $-$ $\hat{\rho}_{h}|\},$
where $\mathfrak{W}_{T,h}^{(k)}$ are possibly stochastic, $\lim
\inf_{T\rightarrow \infty }\min_{\mathcal{H}_{T},\mathcal{K}_{T}}\mathfrak{W}%
_{T,h}^{(k)}$ $>$ $0$ $a.s.$, and $\max_{\mathcal{H}_{T},\mathcal{K}_{T}}|%
\mathfrak{W}_{T,h}^{(k)}$ $-$ $\mathfrak{W}_{h}^{(k)}|$ $\overset{p}{%
\rightarrow }$ $0$ where non-stochastic $\mathfrak{W}_{h}^{(k)}$ satisfy $%
\min_{\mathcal{H}_{T},\mathcal{K}_{T}}\mathfrak{W}_{h}^{(k)}$ $>$ $0$.
Choices include Ljung-Box type weights, or an inverted non-parametric
standard deviation estimator, cf. \cite{HillMotegi2020}.

Consider the latter, and define a sample covariance function $\hat{v}%
_{T}(i;h,k)$ $\equiv $ $1/T\sum_{t=1}^{T-h-i}\hat{z}_{t}(h,k)$ $\times $ $%
\hat{z}_{t+i}(h,k)$ where 
\begin{equation*}
\hat{z}_{t}(h,k)\equiv \left\{ X_{t}X_{t+h}B_{k}(t)-\frac{1}{T}%
\sum_{t=1}^{T-h}X_{t}X_{t+h}B_{k}(t)\right\} .
\end{equation*}%
Under fourth order stationarity under the null, and because $\hat{\gamma}%
_{0} $ only operates as a scale asymptotically, cf. Theorem \ref%
{thm:clt_max_cov_mix_H0} below, the weights are $\mathfrak{W}_{T,h}^{(k)}$ $%
= $ $1/\mathcal{\hat{V}}_{T}(h,k)$ where, e.g., 
\begin{equation}
\mathcal{\hat{V}}_{T}^{2}(h,k)=\hat{\gamma}_{0}^{-2}\left\{ \hat{v}%
_{T}(0;h,k)+2\sum_{i=1}^{T-h-1}\mathcal{K}(i/\beta _{T})\hat{v}%
_{T}(i;h,k)\right\}  \label{V_hac}
\end{equation}%
with symmetric, square integrable kernel function $\mathcal{K}$ $:$ $\mathbb{%
R}$ $\rightarrow $ $[-1,1]$ satisfying $\mathcal{K}(0)$ $=$ $1$,\footnote{%
See, for example, class $\mathfrak{K}_{2}$ in \cite{Andrews1991}, or class $%
\mathfrak{K}$ in \citet[Assumption 1]{deJongDavidson2000}.} and bandwidth $%
\beta _{T}$ $\rightarrow $ $\infty $ where $\beta _{T}$ $=$ $o(T)$.

A penalized weighted version is thus: 
\begin{equation}
\mathcal{M}_{T}^{(w,p)}\equiv \max_{0\leq h\leq \mathcal{H}_{T},1\leq k\leq 
\mathcal{K}_{T}}\left\{ \sqrt{T}\mathfrak{W}_{T,h}^{(k)}\left\vert \hat{\rho}%
_{h}^{(k)}-\hat{\rho}_{h}\right\vert -\mathcal{P}(h,k)\right\} .  \label{Mwp}
\end{equation}%
In Monte Carlo work we study $\mathcal{M}_{T}$, $\mathcal{M}_{T}^{(p)}$,\
and $\mathcal{M}_{T}^{(w,p)}$ with\ Walsh or Haar bases, various penalties,
and/or an inverted standard deviation weight or Ljung-Box weight. We find $%
\mathcal{P}(h,k)$ $=$ $p_{h}$ $+$ $q_{k}$ where $p_{h}$ $=$ $(h$ $+$ $%
1)^{a}/2$ and $q_{k}$ $=$ $k^{a}/2$ with $a$ $=$ $[1/8,1/2]$, or $\mathcal{P}%
(h,k)$ $=$ $\sqrt{(h+1)k}$, promote accurate empirical size but generally
does not lead to dominant power, and may lead to decreased power in some
cases. Conversely, inverted standard error weights $\mathfrak{W}_{T,h}^{(k)}$
generally lead to over-sized tests, and Ljung-Box weights do not offer an
advantage under either hypothesis. Overall a non-penalized and non-weighted
statistic dominates.

Finally, a power improvement may be yielded by combining bases with uniquely
defined systematic samples. Let $\mathcal{M}_{T}(\mathcal{B}_{j})$ be
max-statistics based on $\mathcal{J}$ $\in $ $\mathbb{N}$ orthonormal bases $%
\mathcal{B}_{j,k}(x)$, $j$ $=$ $1,...,\mathcal{J}$. \ Then, ignoring
penalties and weights to ease notation, define a so-called \textquotedblleft 
\textit{max-max-statistic}\textquotedblright :%
\begin{equation}
\mathcal{\check{M}}_{T}\equiv \max_{1\leq j\leq \mathcal{J}}\left\{ \mathcal{%
M}_{T}(\mathcal{B}_{j})\right\} .  \label{M_combo}
\end{equation}%
We study $\mathcal{\check{M}}_{T}\equiv \max \left\{ \mathcal{M}_{T}(%
\mathcal{W}),\mathcal{M}_{T}(\Psi )\right\} $ in simulation work, where $%
\mathcal{M}_{T}(\mathcal{W})$ and $\mathcal{M}_{T}(\Psi )$ use Walsh and
composite Haar bases. An asymptotic theory for $\mathcal{\check{M}}_{T}$ and
its bootstrapped p-value follow directly from results given below and the
mapping theorem since $\mathcal{J}$\ is a finite constant. Other basis
combinations are clearly feasible. Consider discretized bases $B_{j,k}(t)$
and the set $\{\bar{B}_{\bar{k}}(t)\}_{\bar{k}=1}^{\mathcal{\bar{K}}}$ $=$ $%
\{B_{j,k}(t)$ $:$ $j$ $\in $ $\mathcal{J}^{\ast };k$ $\in $ $\mathcal{K}%
^{\ast }\}$ where $\mathcal{J}^{\ast }$ and $\mathcal{K}^{\ast }$ are index
subsets of $\{1,...,\mathcal{J}\}$ and $\{1,...,\mathcal{K}\}$ yielding
unique $\mathcal{B}_{j,k}(x)$ $\forall x$. Test statistics can then be
derived from $\{\bar{B}_{\bar{k}}(t)\}_{\bar{k}=1}^{\mathcal{\bar{K}}}$.

\section{Asymptotic theory\label{sec:max_corr_asym}}

Write 
\begin{eqnarray}
&&z_{t}(h,k)\equiv X_{t}X_{t+h}B_{k}(t)-E\left[ X_{t}X_{t+h}\right] B_{k}(t)
\label{zthk} \\
&&\mathcal{Z}_{T}(h,k)\equiv \frac{1}{\sqrt{T}}\sum_{t=1}^{T-h}z_{t}(h,k),
\label{ZT}
\end{eqnarray}%
and define a variance function $\sigma _{T}^{2}(h,k)$ $\equiv $ $E[\mathcal{Z%
}_{T}^{2}(h,k)]$. In the general case $E[X_{t}]$ $\in $ $\mathbb{R}$ replace 
$X_{t}$ with $X_{t}$ $-$ $E[X_{t}]$.

The main contributions of this section deliver a class of sequences $\{%
\mathcal{H}_{T},\mathcal{K}_{T}\}$, and an array of random variables $\{%
\boldsymbol{Z}_{T}(h,k)$ $:$ $T$ $\in $ $\mathbb{N}\}_{h\geq 0,k\geq 1}$
normally distributed $\boldsymbol{Z}_{T}(h,k)$ $\sim $ $N(0,\sigma
_{T}^{2}(h,k))$, such that the Kolmogorov distance 
\begin{equation}
\rho _{T}\equiv \sup_{z\geq 0}\left\vert P\left( \max_{\mathcal{H}_{T},%
\mathcal{K}_{T}}\left\vert \mathcal{Z}_{T}(h,k)\right\vert \leq z\right)
-P\left( \max_{\mathcal{H}_{T},\mathcal{K}_{T}}\left\vert \boldsymbol{Z}%
_{T}(h,k)\right\vert \leq z\right) \right\vert \rightarrow 0.
\label{Gauss_apprx}
\end{equation}%
We work under mixing or physical dependence settings. The approximation does
not require standardized $\mathcal{Z}_{T}$ and $\boldsymbol{Z}_{T}$ in view
of non-degeneracy Assumption \ref{assum_1abcd}.c below. We then apply the
approximation to the max-correlation difference statistic. 

\subsection{Mixing}

Define $\sigma $-fields $\mathcal{F}_{T,t}^{\infty }$ $\equiv $ $\sigma
\left( X_{\tau }:\tau \geq t\right) $ and $\mathcal{F}_{T,-\infty }^{t}$ $%
\equiv $ $\sigma \left( X_{\tau }:\tau \leq t\right) ,$ and $\alpha $-mixing
coefficients \citep{Rosenblatt1956}, $\alpha (l)$ $\equiv $ $\sup_{t\in 
\mathbb{Z}}\sup_{\mathcal{A}\subset \mathcal{F}_{T,-\infty }^{t},\mathcal{B}%
\subset \mathcal{F}_{T,t+l}^{\infty }}|\mathcal{P}(\mathcal{A}\cap \mathcal{B%
})$ $-$ $\mathcal{P}\left( \mathcal{A}\right) \mathcal{P}\left( \mathcal{B}%
\right) |$, for $l$ $>$ $0.$ We work in the setting of \cite%
{ChangJiangShao2023}, cf. \cite{ChangChenWu2024}, who deliver high
dimensional central limit theorems for possibly non-stationary mixing
sequences or under a physical dependence setting similar to \cite%
{ZhangWu2017}. \citet[Appendix
B]{Chernozhukov_etal_manymom_2014}, cf. 
\citet[Supplemental
Appendix]{Chernozhukov_etal2019}, allow for \textit{almost surely} bounded
stationary $\beta $-mixing processes, while \cite{ZhangWu2017} extend
results in \cite{Chernozhukov_etal2013} to a large class of dependent
stationary processes. Stationarity is not suitable here because even under
the null we want to allow for global non-stationarity, and boundedness is
typically too restrictive for many financial and macroeconomic time series.

\begin{assumption}
\label{assum_1abcd} $\ \ \medskip $\newline
$a$. (\emph{weak dependence}): $\alpha (l)$ $\leq $ $K_{1}\exp
\{-K_{2}l^{\phi }\}$ for some universal constants $\phi ,K_{1},K_{2}$ $>$ $0$%
.$\medskip $\newline
$b.$ (\emph{subexponential tails}): $\max_{1\leq t\leq T}P(|X_{t}|$ $>$ $c)$ 
$\leq $ $\vartheta _{1}\exp \{-\vartheta _{2}c^{\varpi }\}$ for some
universal constants $\varpi ,\vartheta _{1},\vartheta _{2}$ $>$ $0$.$%
\medskip $\newline
$c.$ (\emph{nondegeneracy}): $\lim \inf_{T\rightarrow \infty }E[\mathcal{Z}%
_{T}^{2}(h,k)]$ $>$ $0$ $\forall (h,k)$.$\medskip $\newline
$d.$ (\emph{orthonormal basis}): $\{\mathcal{B}_{k}(x)$ $:$ $0$ $\leq $ $k$ $%
\leq $ $\mathcal{K}\}$ forms a complete orthonormal basis on $\mathcal{L}%
[0,1)$; $\mathcal{B}_{k}(x)$ $\in $ $\{-1,1\}$ on $[0,1)$; $%
|\sum_{t=1}^{T}B_{k}(t)|$ $=$ $O(\eta (k))$ for some positive strictly
monotonic function $\eta $ $:$ $\mathbb{R}_{+}$ $\rightarrow $ $\mathbb{R}%
_{+}$, $\eta (k)$ $\nearrow $ $\infty $ as $k$ $\rightarrow $ $\infty $.
\end{assumption}

\begin{remark}
\label{rm:z}\normalfont A version of ($a$)-($c$) are imposed in %
\citet[Conditions AS1-AS3]{ChangJiangShao2023}. Under ($a$) their high
dimensional mixing condition AS2 applies to $[z_{t}(h,k)]_{h=1,k=0}^{%
\mathcal{H}_{T},\mathcal{K}_{T}}$. ($b$) trivially generalizes their AS1,
which sets $\vartheta _{1}$ $\geq $ $1$ and $\varpi $ $\in $ $(0,1]$ for
notational convenience. It covers traditional sub-exponential tails, and
slower decay such that a moment generating function may not exist 
\citep[Proposition
2.7.1]{Vershynin2018}, while permitting general forms of global
nonstationarity. ($c$) yields AS3.
\end{remark}

\begin{remark}
\normalfont Nondegeneracy ($c$) is common in the time series literature %
\citep[e.g.][Theorem 1]{Doukhan1994}, in particular for non-standardized
statistics involving nonstationary sequences.\footnote{%
See, e.g., \citet[Assumption
2.a]{deJong1997}, but also see 
\citet[Th\underline{}eorem
19.1]{Billingsley1999}.}. It classically rules out degenerate dispersion and
deviant negative co-dependence within the sequence $\{X_{t}X_{t+h}$ $-$ $%
E[X_{t}X_{t+h}]\}_{t=1}^{T-h}$. Simply note $E[\mathcal{Z}_{T}^{2}(h,k)]$ $=$
$((T$ $-$ $h)/T)$ $\times $ $E[(\boldsymbol{\lambda }_{T-h}^{\prime }%
\boldsymbol{\check{X}}_{T-h})^{2}]$ where $\boldsymbol{\check{X}}_{T-h}$ $%
\equiv $ $[X_{1}X_{1+h}$ $-$ $E[X_{1}X_{1+h}]]_{t=1}^{T-h}$ and $\boldsymbol{%
\lambda }_{T-h}$ $\equiv $ $(T$ $-$ $h)^{-1/2}[B_{k}(1),...,B_{k}(T$ $-$ $%
h)] $. Notice $\boldsymbol{\lambda }_{T-h}^{\prime }\boldsymbol{\lambda }%
_{T-h}$ $=$ $1$ since $B_{k}^{2}(t)$ $=$ $1$. Thus ($c$) is satisfied by a
classic positive definiteness property: $\inf_{\lambda ^{\prime }\lambda
=1}E[(\lambda ^{\prime }\boldsymbol{\check{X}}_{T-h})^{2}]$ $>$ $0$ $\forall
(h,k)$ and $\forall T$ $\geq $ $\underline{T}$ and some $\underline{T}$ $\in 
$ $\mathbb{N}$. For example, impose fourth order stationarity (and therefore
the null), and white noise $E[X_{t}X_{t+h}]$ $=$ $0$ $\forall h$ $\geq $ $1$
to reduce notation. Now define fourth order correlation coefficients $%
r(a,b,c,d)$ $\equiv $ $E\left[ X_{a}X_{b}X_{c}X_{d}\right]
/E[X_{a}^{2}X_{c}^{2}]$. Then by expanding $E[(\lambda ^{\prime }\boldsymbol{%
\check{X}}_{T-h})^{2}]$, ($c$) holds under pointwise non-degeneracy $%
E[X_{1}^{2}X_{1+h}^{2}]$ $>$ $0$, and $\inf_{\lambda ^{\prime }\lambda
=1}\{1 $ $+$ $\sum_{i=1}^{T-h-1}r(0,h,i,i$ $+$ $h)\sum_{t=1}^{T-h-i}\lambda
_{t}\lambda _{t+i}\}$ $>$ $0$ $\forall (h,k)$, $\forall T$ $\geq $ $%
\underline{T}$, ruling out deviant negative linear dependence. See also the
discussion in \citet[p.
4-5]{ChangChenWu2024}. We cannot, however, impose fourth order stationarity 
\emph{broadly}, and thus the preceding sufficient conditions, because that
rules out an asymptotic analysis under local or global alternatives, cf.
Section \ref{sec:max_corr_H1}.
\end{remark}

\begin{remark}
\normalfont Assumption \ref{assum_1abcd} reveals a trade-off vis-\`{a}-vis
JWW. We allow for nonlinear $\{X_{t}\}$ with possibly non-iid errors, and
possibly global nonstationarity under the null, but $X_{t}$ must have
exponentially decaying tails and geometric dependence. The former rules out
conventional GARCH processes (which lack higher moments), but includes
GARCH-type processes with errors that have bounded support. JWW focus
exclusively on linear processes $X_{t}=\sum_{i=0}^{\infty }\psi _{i}Z_{t-i}$
with iid $Z_{t}$ where $E|Z_{t}|^{4\upsilon }$ $<$ $\infty $ for some $%
\upsilon $ $>$ $1$, excluding important nonlinear and conditionally
heteroscedastic processes. They impose $\psi _{i}=O(1/[i(\ln i)^{1+\kappa
}]) $ for some $\kappa $ $>$ $0$ and strict stationarity under the null,
yielding $\sum_{h=1}^{\infty }|\gamma _{h}|$ $<$ $\infty $. Thus JWW allow
for hyperbolic \emph{and} geometric memory decay and the possible
nonexistence of higher moments. Under physical dependence Assumption \ref%
{assum_1abcd}.a$^{\ast }$ below, however, we allow for hyperbolic memory.
\end{remark}

\begin{remark}
\normalfont The bound $|\sum_{t=1}^{T}B_{k}(t)|$ $=$ $O(\eta (k))$ in ($d$)\
is generally driven by the number of zero crossings on $[0,1)$ in the
underlying smooth basis function $\mathcal{B}_{k}(x)$. Indeed, by Lemma 3 in
JWW, Walsh $\mathcal{W}_{k}(t)$ exhibit up to $k$ zero crossings, and $%
|\sum_{t=1}^{T}\mathcal{W}_{k}(t)|$ $\leq $ $k$ $+$ $1$ hence $\eta (k)$ $=$ 
$k$. Conversely, Haar composite $\Psi _{k}(t)$ exhibits up to $2^{k}$ zero
crossings, and $|\sum_{t=1}^{T}\Psi _{k}(t)|$ $=$ $O(2^{k})$ by Lemma B.1 in 
\cite{supp_mat_covstat_2024}, hence $\eta (k)$ $=$ $2^{k}$.
\end{remark}

\begin{lemma}
\label{lm:clt_max_cov_mix}Under Assumption \ref{assum_1abcd} we have $\rho
_{T}$ $\lesssim $ $\mathcal{H}_{T}^{1/2}\left( \ln \left( \mathcal{H}%
_{T}\right) \right) ^{7/6}/T^{1/9}$ $\rightarrow $ $0$ for any sequences $\{%
\mathcal{H}_{T},\mathcal{K}_{T}\}$ with $0$ $\leq $ $\mathcal{H}_{T}$ $\leq $
$T$ $-$ $1$, $\mathcal{H}_{T}$ $=$ $O(T^{1/9}(\ln (T))^{1/3})$, $\mathcal{K}%
_{T}$ $=$ $o(T^{\kappa })$ for some finite $\kappa $ $>$ $0$, and $\eta (%
\mathcal{K}_{T})$ $=$ $o(\sqrt{T})$ where $\eta (\cdot )$ is the Assumption %
\ref{assum_1abcd}.d discrete basis summand bound. In this case $\max_{%
\mathcal{H}_{T},\mathcal{K}_{T}}|\mathcal{Z}_{T}(h,k)|$ $\overset{d}{%
\rightarrow }$ $\max_{h,k\in \mathbb{N}}|\boldsymbol{Z}(h,k)|$ where $%
\boldsymbol{Z}(h,k)$ $\sim $ $N(0,\lim_{T\rightarrow \infty }\sigma
_{T}^{2}(h,k))$ and $\lim_{T\rightarrow \infty }\sigma _{T}^{2}(h,k)$ $<$ $%
\infty $.
\end{lemma}

\begin{remark}
\label{rm:KT}\normalfont We require the orthonormal basis $\mathcal{B}(x)$
bound function $\eta (\cdot )$ and maximum systematic sample counter $%
\mathcal{K}_{T}$ to satisfy $\eta (\mathcal{K}_{T})$ $=$ $o(\sqrt{T})$ to
ensure the mean summation $\mathcal{S}_{T}^{(k)}(h)$ $\equiv $ $1/\sqrt{T}%
\sum_{t=1}^{T-h}E[X_{t}X_{t+h}]B_{k}(t)$ is negligible in the proof of
Theorem \ref{thm:clt_max_cov_mix_H0} below. Simply note that under $H_{0}$
and Assumption \ref{assum_1abcd}.d, $|\mathcal{S}_{T}^{(k)}(h)|$ $\leq $ $%
\gamma _{h}|1/\sqrt{T}\sum_{t=1}^{T-h}B_{k}(t)|$ $\leq $ $\gamma _{h}\eta
(k)/\sqrt{T}$. Thus $\max_{1\leq k\leq \mathcal{K}_{T}}|\mathcal{S}%
_{T}^{(k)}|$ $\leq $ $\gamma _{h}\eta (\mathcal{K}_{T})/\sqrt{T}$ $%
\rightarrow $ $0$ when $\eta (\mathcal{K}_{T})=o(\sqrt{T})$. In a time
series setting $\mathcal{H}_{T}$ $=$ $o(T)$ must hold to ensure consistency
of sample autocovariances. Specifically $\mathcal{H}_{T}$ $=$ $O(T^{1/9}(\ln
(T))^{1/3})$ and $\mathcal{K}_{T}$ $=$ $o(T^{\kappa })$ for any $\kappa $ $>$
$0$ yield $\rho _{T}$ $\lesssim $ $\mathcal{H}_{T}^{1/2}\left( \ln \left( 
\mathcal{H}_{T}\right) \right) ^{7/6}/T^{1/9}$. $\mathcal{K}_{T}$ $=$ $%
o(T^{\kappa })$ is implied by $\eta (\mathcal{K}_{T})=o(\sqrt{T})$ for Walsh
and Haar functions: see below. $\mathcal{H}_{T}$ $=$ $O(T^{1/9}(\ln
(T))^{1/3})$ is necessary due to high dimensional lagging: $%
[z_{t}(h,k)]_{h=1,k=0}^{\mathcal{H}_{T},\mathcal{K}_{T}}$ is $\sigma
(X_{\tau }$ $:$ $\tau $ $\leq $ $t$ $+$ $\mathcal{H}_{T})$-measurable and
therefore has mixing coefficients $\alpha (|l-\mathcal{H}_{T}|_{+})$,
greatly impacting feasible $\{\mathcal{H}_{T}\}$ 
\citep[cf.][Proposition
3]{ChangJiangShao2023}.
\end{remark}

\begin{remark}
\normalfont Walsh functions $\mathcal{W}_{k}(t)$ have $\eta (k)$ $=$ $k$
hence $\mathcal{K}_{T}$ $=$ $o(\sqrt{T})$, while Haar composite $\Psi
_{k}(t) $ have $\eta (k)$ $=$ $2^{k}$ hence $\mathcal{K}_{T}$ $=$ $o(\ln
(T)) $, yielding $\mathcal{K}_{T}$ $=$ $o(T^{\kappa })$ respectively for
some, or any, $\kappa $ $>$ $0$.
\end{remark}

Now define $\sigma ^{2}(h,k)$ $\equiv $ $\lim_{T\rightarrow \infty }\sigma
_{T}^{2}(h,k).$ Under $H_{0}$ and Assumption \ref{assum_1abcd}, $\sigma
^{2}(h,k)$ $\in $ $(0,\infty )$. We now have a limit theory for the
max-correlation difference.

\begin{theorem}
\label{thm:clt_max_cov_mix_H0}Let $H_{0}$ and Assumption \ref{assum_1abcd}\
hold, and let $\mathcal{H}_{T},\mathcal{K}_{T}$ $\rightarrow $ $\infty $.
Let $\{\boldsymbol{Z}(h,k)$ $:$ $h,k$ $\in $ $\mathbb{N}\}$ be a zero mean
Gaussian process with $\boldsymbol{Z}(h,k)$ $\sim $ $N(0,\sigma ^{2}(h,k))$.
Then it holds that $\mathcal{M}_{T}$ $\overset{d}{\rightarrow }$ $\gamma
_{0}^{-1}\max_{h,k\in \mathbb{N}}\left\vert \boldsymbol{Z}(h,k)\right\vert $
for any $\{\mathcal{H}_{T},\mathcal{K}_{T}\}$ with $0$ $\leq $ $\mathcal{H}%
_{T}$ $\leq $ $T$ $-$ $1$, $\mathcal{H}_{T}$ $=$ $O(T^{1/9}(\ln (T))^{1/3})$%
, $\mathcal{K}_{T}$ $=$ $o(T^{\kappa })$ for some finite $\kappa $ $>$ $0$,
and $\eta (\mathcal{K}_{T})$ $=$ $o(\sqrt{T})$.
\end{theorem}

\begin{remark}
\normalfont Consider the weighted/penalized version $\mathcal{M}_{T}^{(w,p)}$
in (\ref{Mwp}), and assume the weights satisfy $\lim \inf_{T\rightarrow
\infty }\inf_{\mathcal{H}_{T},\mathcal{K}_{T}}\mathfrak{W}_{T,h}^{(k)}$ $>$ $%
0$ $a.s.$, and $\max_{\mathcal{H}_{T},\mathcal{K}_{T}}|\mathfrak{W}%
_{T,h}^{(k)}$ $-$ $\mathfrak{W}_{h}^{(k)}|$ $\overset{p}{\rightarrow }$ $0$
where non-stochastic $\mathfrak{W}_{h}^{(k)}$ satisfy $\inf_{h,k\in \mathbb{N%
}}\mathfrak{W}_{h}^{(k)}$ $>$ $0$. The penalty functions $(p_{v},q_{v})$ are
positive, monotonically increasing and bounded on compact sets. Then from
arguments used to prove Theorem \ref{thm:clt_max_cov_mix_H0}, it follows:%
\begin{equation*}
\mathcal{M}_{T}^{(w,p)}\overset{d}{\rightarrow }\gamma _{0}^{-1}\max_{k\in 
\mathbb{N}}\left[ \max_{h\in \mathbb{N}}\left\{ \mathfrak{W}%
_{h}^{(k)}\left\vert \boldsymbol{Z}(h,k)\right\vert -p_{h}\right\} -q_{k}%
\right]
\end{equation*}

Now suppose we standardize with $\mathfrak{W}_{T,h}^{(k)}$ $=$ $1/\mathcal{%
\hat{V}}_{T}(h,k)$ with HAC estimator $\mathcal{\hat{V}}_{T}^{2}(h,k)$ in (%
\ref{V_hac}), and kernel function $\mathcal{K}(\cdot )$ belonging to class $%
\mathfrak{K}$ in \citet[Assumption 1]{deJongDavidson2000}, or class $%
\mathfrak{K}_{2}$ $\supset $ $\mathfrak{K}$ in \cite{Andrews1991}. \cite%
{deJongDavidson2000} allow for possibly globally nonstationary mixing
sequences (or non-mixing satisfying a \emph{near epoch dependence}
property). In their environment with bandwidth $\beta _{T}$ $=$ $o(T)$ we
have $\mathfrak{W}_{T,h}^{(k)}$ $>$ $0$ $a.s.$ and $\mathfrak{W}_{T,h}^{(k)}$
$\overset{p}{\rightarrow }$ $\mathfrak{W}_{h}^{(k)}$ $=$ $1/\mathcal{V}(h,k)$
where $\mathcal{V}^{2}(h,k)$ $=$ $\gamma _{0}^{-1}\lim_{T\rightarrow \infty
}\sigma _{T}^{2}(h,k)$. Uniformity $\max_{\mathcal{H}_{T},\mathcal{K}_{T}}|%
\mathfrak{W}_{T,h}^{(k)}$ $-$ $\mathfrak{W}_{h}^{(k)}|$ $\overset{p}{%
\rightarrow }$ $0$ can be proved using theory developed in this section, and
Section \ref{sec:bootstrap}, omitted here for space considerations.
\end{remark}

\subsection{Physical dependence}

Now consider the setting of \cite{Wu2005}, and as in \cite{ChangChenWu2024}
assume $X_{t}=g_{t}(\epsilon _{t},\epsilon _{t-1},\ldots )$ for measurable $%
g_{t}(\cdot )$ that may depend on $t$, and iid stochastic $\{\epsilon
_{t}\}_{t\in \mathbb{Z}}$. Let $\{\epsilon _{t}^{\prime }\}_{t\in \mathbb{Z}%
} $ be an independent copy of $\{\epsilon _{t}\}_{t\in \mathbb{Z}}$, and let 
$X_{t}^{\prime }(m)$ $=$ $g_{t}(\epsilon _{t},\epsilon _{t-1},\ldots
,\epsilon _{t-m}^{\prime },\epsilon _{t-m-1},...)$ be coupled versions, $m$ $%
=$ $0,1,2,..$. Define $\theta _{t}^{(p)}(m)$ $\equiv $ $||X_{t}$ $-$ $%
X_{t}^{\prime }(m)||_{p}$. We say $X_{t}$ is $\mathcal{L}_{p}$-\textit{%
physical dependent} when $\lim_{m\rightarrow \infty }\lim \sup_{T\rightarrow
\infty }\max_{1\leq t\leq T}\theta _{t}^{(p)}(m)$ $=$ $0$. We make the
following assumption for ease of exposition.\medskip \newline
\textbf{Assumption \ref{assum_1abcd}.a}$^{\ast }$\textbf{.} (\emph{weak
dependence}): \textit{$X_{t}$ is $\mathcal{L}_{p}$-physical dependent for
some $p$ $\geq $ $8$, with $\theta _{t}^{(p)}(m)$ $\leq $ $d_{t}^{(p)}\psi
_{m}$ where $\psi _{m}$ $=$ $O(m^{-\lambda -\iota })$ for some }size\textit{%
\ $\lambda $ $\geq $ $1$.}

\begin{remark}
\normalfont The bound is similar to mixingale and near epoch dependence
constructs: $d_{t}^{(p)}$ captures time-dependent heterogeneity, and $\psi
_{m}$ dictates memory decay %
\citep[e.g.][]{McLeish1975,DedeckerDoukhan_etal2007}. Such a representation
holds for many linear and nonlinear processes \citep[e.g.][]{Wu2005}. (a$%
^{\ast }$) allows for hyperbolic memory, contrary to mixing (a).
\end{remark}

\begin{theorem}
\label{thm:clt_max_cov_phys_H0}Let $H_{0}$ and Assumption \ref{assum_1abcd}.a%
$^{\ast }$,b,c,d\ hold, and let $\mathcal{H}_{T},\mathcal{K}_{T}$ $%
\rightarrow $ $\infty $. Let $\{\boldsymbol{Z}(h,k)$ $:$ $h,k$ $\in $ $%
\mathbb{N}\}$ be a zero mean Gaussian process with $\boldsymbol{Z}(h,k)$ $%
\sim $ $N(0,\sigma ^{2}(h,k))$. Then it holds that $\mathcal{M}_{T}$ $%
\overset{d}{\rightarrow }$ $\gamma _{0}^{-1}\max_{h,k\in \mathbb{N}%
}\left\vert \boldsymbol{Z}(h,k)\right\vert $ for any $\{\mathcal{H}_{T},%
\mathcal{K}_{T}\}$ with $0$ $\leq $ $\mathcal{H}_{T}$ $\leq $ $T$ $-$ $1$, $%
\mathcal{H}_{T}$ $=$ $o(T)$, $\mathcal{K}_{T}$ $=$ $o(T^{\kappa })$ for some
finite $\kappa $ $>$ $0$, and $\eta (\mathcal{K}_{T})$ $=$ $o(\sqrt{T})$.
\end{theorem}

\begin{remark}
\normalfont The result is based on a high dimensional central limit theorem
in \citet[Theorem 3]{ChangChenWu2024}, requiring physical dependence for $%
z_{t}(h,k)$ \emph{uniformly} over $(h,k)$. Under \emph{joint} mixing for $%
[z_{t}(h,k)]_{h=0,k=1}^{\mathcal{H}_{T},\mathcal{K}_{T}}$ in 
\citet[Proposition
3]{ChangJiangShao2023}, however, we only achieve $\mathcal{H}_{T}$ $=$ $%
o(T^{1/9})$ because $[z_{t}(h,k)]_{h=0,k=1}^{\mathcal{H}_{T},\mathcal{K}%
_{T}} $ has coefficients $\alpha (|l-\mathcal{H}_{T}|_{+})$.
\end{remark}

\subsection{Max-correlation difference under $H_{1}$ \label{sec:max_corr_H1}}

The correlation difference expands to:%
\begin{equation}
\sqrt{T}(\hat{\rho}_{h}^{(k)}-\hat{\rho}_{h})=\frac{1}{\hat{\gamma}_{0}}%
\frac{1}{\sqrt{T}}\sum_{t=1}^{T-h}\left( X_{t}X_{t+h}-E\left[ X_{t}X_{t+h}%
\right] \right) B_{k}(t)+\frac{1}{\hat{\gamma}_{0}}\frac{1}{\sqrt{T}}%
\sum_{t=1}^{T-h}E\left[ X_{t}X_{t+h}\right] B_{k}(t).  \label{corrdiff}
\end{equation}%
Under either hypothesis $1/\sqrt{T}\sum_{t=1}^{T-h}(X_{t}X_{t+h}$ $-$ $E%
\left[ X_{t}X_{t+h}\right] )B_{k}(t)$ is asymptotically normal. For the
sample variance, we similarly have under either hypothesis and Assumption %
\ref{assum_1abcd}:%
\begin{equation*}
\sqrt{T}\left( \hat{\gamma}_{0}-\frac{1}{T}\sum_{t=1}^{T}E\left[ X_{t}^{2}%
\right] \right) =\frac{1}{\sqrt{T}}\sum_{t=1}^{T}\left( X_{t}^{2}-E\left[
X_{t}^{2}\right] \right) =O_{p}(1).
\end{equation*}%
Hence, $\hat{\gamma}_{0}$ $=$ $g_{0}$ $+$ $O_{p}(1/\sqrt{T})$ assuming
existence of $g_{0}$ $\equiv $ $\lim_{T\rightarrow \infty
}1/T\sum_{t=1}^{T}E[[X_{t}^{2}]$. See below for derivations of $g_{0}$ under
local and global alternatives.

In order to handle $1/\sqrt{T}\sum_{t=1}^{T-h}E[X_{t}X_{t+h}]B_{k}(t)$ in (%
\ref{corrdiff}), we need a representation of a non-stationary covariance for
fixed and local alternatives. Let $\gamma _{h}(u)$ be the time varying
autocovariance function on $[0,1]$. In the framework of locally stationary
processes \citep[cf.][]{Dahlhaus1997,Dahlhaus2009}, we may state the global
alternative hypothesis as%
\begin{equation}
H_{1}:\int_{0}^{1}\left( \gamma _{h}(u)-\int_{0}^{1}\gamma _{h}(v)dv\right)
^{2}du>0\text{ for some }h\geq 0.  \label{H1_gam(u)}
\end{equation}%
Thus under $H_{1}$ there exists a lag $h$ and subset $\mathcal{S}_{h}$ $%
\mathcal{\subset }$ $[0,1]$ with positive Lebesgue measure such that $\gamma
_{h}(u)$ $\neq $ $\int_{0}^{1}\gamma _{h}(v)dv$ on $\mathcal{S}_{h}$; hence $%
\gamma _{h}(u)$ is not \textit{almost everywhere} constant on $[0,1]$.

Now, by completeness of $\{\mathcal{B}_{k}(u)$ $:$ $0$ $\leq $ $k$ $\leq $ $%
\mathcal{K}\}$ under Assumption \ref{assum_1abcd}.d, we may write $\gamma
_{h}(u)$ $=$ $\sum_{k=0}^{\infty }\omega _{h,k}\mathcal{B}_{k}(u)$ $=$ $%
\omega _{h,0}$ $+$ $\sum_{k=1}^{\infty }\omega _{h,k}\mathcal{B}_{k}(u)$,
where $\omega _{h,k}$ $=$ $\int_{0}^{1}\gamma _{h}(u)\mathcal{B}_{k}(u)$ by
orthonormality. Hence, under $H_{1}$ and orthonormality, for some $h\geq 0.$ 
\begin{equation*}
\int_{0}^{1}\left( \gamma _{h}(u)-\int_{0}^{1}\gamma _{h}(v)dv\right)
^{2}du=\int_{0}^{1}\left( \sum_{k=1}^{\infty }\omega _{h,k}\mathcal{B}%
_{k}(u)\right) ^{2}du=\sum_{k=1}^{\infty }\omega _{h,k}^{2}>0,
\end{equation*}%
which yields $\max_{h,k\in \mathbb{N}}|\int_{0}^{1}\gamma _{h}(u)\mathcal{B}%
_{k}(u)|$ $>$ $0$ under $H_{1}$.

A sequence of local alternatives with $\sqrt{T}$-drift logically follows:%
\begin{equation}
H_{1}^{L}:E\left[ X_{t}X_{t+h}\right] =\gamma _{h}+c_{h}(t/T)/\sqrt{T},\text{
}  \label{H1L}
\end{equation}%
where $\gamma _{h}$ is a constant for each $h$, $\max_{h\in \mathbb{N}%
}|\gamma _{h}|$ $\leq $ $K$ $<$ $\infty $, and $c_{h}$ $:$ $[0,1]$ $%
\rightarrow $ $\mathbb{R}$ are integrable functions on $[0,1]$ uniformly\
over $h$ (i.e. $\sup_{h\in \mathbb{N}}|\int_{0}^{1}c_{h}(u)du)|$ $<$ $\infty 
$), that satisfy (\ref{H1_gam(u)}). Thus, under local alternative (\ref{H1L}%
), by the preceding discussion:%
\begin{equation}
\liminf_{T\rightarrow \infty }\max_{h,k\in \mathbb{N}}\left\vert
\int_{0}^{1}c_{h}\left( u\right) \mathcal{B}_{k}(u)du\right\vert >0\text{.}
\label{cW}
\end{equation}%
In order to ensure $\min_{t\in \mathbb{Z}}E[X_{t}^{2}]$ $>$ $0$, assume $%
\gamma _{0}$ $>$ $0$ and $c_{0}(u)$ $\geq $ $0$ \textit{almost everywhere}.
Notice $\lim_{T\rightarrow \infty }|T^{-1}\sum_{t=1}^{T}c_{0}(t/T)|$ $=$ $%
|\int_{0}^{1}c_{0}(u)du|$ $<$ $\infty $ yields under $H_{1}^{L}$: 
\begin{equation*}
g_{0}\equiv \lim_{T\rightarrow \infty }\frac{1}{T}%
\sum_{t=1}^{T}E[X_{t}^{2}]=\gamma _{0}+\lim_{T\rightarrow \infty }\frac{1}{%
\sqrt{T}}\frac{1}{T}\sum_{t=1}^{T}c_{0}(t/T)=\gamma _{0}.
\end{equation*}

Under Assumption \ref{assum_1abcd}.d $|\sum_{t=1}^{T}B_{k}(t)|$ $=$ $O(\eta
(k))$, and $\max_{h\in \mathbb{N}}|\gamma _{h}|$ $\leq $ $K$ and $\eta (%
\mathcal{K}_{T})$ $=$ $o(\sqrt{T})$\ by supposition. Hence $\max_{\mathcal{H}%
_{T},\mathcal{K}_{T}}$ $|1/\sqrt{T}\sum_{t=1}^{T-h}B_{k}(t)|$ $=$ $o(1)$.
Thus under $H_{1}^{L}$: 
\begin{eqnarray}
\frac{1}{\sqrt{T}}\sum_{t=1}^{T-h}E\left[ X_{t}X_{t+h}\right] B_{k}(t)
&=&\gamma _{h}\frac{1}{\sqrt{T}}\sum_{t=1}^{T-h}B_{k}(t)+\frac{1}{T}%
\sum_{t=1}^{T-h}c_{h}(t/T)B_{k}(t)  \label{EXXW_HL1} \\
&=&o\left( 1\right) +\frac{1}{T}\sum_{t=1}^{T-h}c_{h}(t/T)B_{k}(t)%
\rightarrow \int_{0}^{1}c_{h}\left( u\right) \mathcal{B}_{k}(u)du,  \notag
\end{eqnarray}%
where here and below $o(1)$, and all subsequent $O_{p}(\cdot )$ and $%
o_{p}(\cdot )$ terms, do not depend on $(h,k)$.

Asymptotics rest on a uniform limit theory over $(h,k)$, which here needs to
extend to the limit in (\ref{EXXW_HL1}). We therefore enhance local
alternative (\ref{H1L}) by assuming $c_{h}(\cdot )$ satisfies for any $\{%
\mathcal{H}_{T},\mathcal{K}_{T}\}$: 
\begin{equation}
\max_{\mathcal{H}_{T},\mathcal{K}_{T}}\left\vert \frac{1}{T}%
\sum_{t=1}^{T-h}c_{h}(t/T)B_{k}(t)-\int_{0}^{1}c_{h}\left( u\right) \mathcal{%
B}_{k}(u)du\right\vert \rightarrow 0.  \label{c(t/T)-c(u)}
\end{equation}

Now define: 
\begin{equation}
\mathcal{C}(h,k)=\int_{0}^{1}c_{h}\left( u\right) \mathcal{B}_{k}(u)du\text{.%
}  \label{C(hk)}
\end{equation}%
Then under $H_{1}^{L}$, $\lim \inf_{T\rightarrow \infty }\max_{h,k\in 
\mathbb{N}}|\mathcal{C}(h,k)|$ $>$ $0$ in view of (\ref{cW}). Use arguments
in the proof of Theorem \ref{thm:clt_max_cov_mix_H0} to yield under $%
H_{1}^{L}$: 
\begin{eqnarray*}
\sqrt{T}\left( \hat{\rho}_{h}^{(k)}-\hat{\rho}_{h}\right)  &=&\frac{1}{g_{0}}%
\frac{1}{\sqrt{T}}\sum_{t=1}^{T-h}\left( X_{t}X_{t+h}-E\left[ X_{t}X_{t+h}%
\right] \right) B_{k}(t) \\
&&+\frac{1}{g_{0}}\left( \int_{0}^{1}c_{h}\left( u\right) \mathcal{B}%
_{k}(u)du+o\left( 1\right) \right) +O_{p}(1/\sqrt{T}).
\end{eqnarray*}%
Hence, by Lemma \ref{lm:clt_max_cov_mix} for any $\{\mathcal{H}_{T},\mathcal{%
K}_{T}\}$ with $\mathcal{H}_{T}$ $=$ $O(T^{1/9}(\ln (T))^{1/3})$, $\mathcal{K%
}_{T}$ $=$ $o(T^{\kappa })$ for some finite $\kappa $ $>$ $0$, and $\eta (%
\mathcal{K}_{T})$ $=$ $o(\sqrt{T})$:%
\begin{equation}
\max_{\mathcal{H}_{T},\mathcal{K}_{T}}\left\vert \sqrt{T}\left( \hat{\rho}%
_{h}^{(k)}-\hat{\rho}_{h}\right) \right\vert \overset{d}{\rightarrow }\frac{1%
}{g_{0}}\max_{h,k\in \mathbb{N}}\left\vert \boldsymbol{Z}(h,k)+\mathcal{C}%
(h,k)\right\vert .  \label{max_p_H1}
\end{equation}%
Thus, the proposed test has non-negligible power under the sequence of $%
\sqrt{T}$-local alternatives (\ref{H1L}) when $c_{h}(\cdot )$ satisfy (\ref%
{H1_gam(u)}), for any complete orthonormal basis in view of (\ref{cW}).
Notice under $H_{0}$ we have $g_{0}$ $=$ $\gamma _{0}$, and $c_{h}(u)$ $=$ $0
$ $\forall u,h$ so that $\mathcal{C}(h,k)$ $=$ $0$ $\forall h,k$, yielding
Theorem \ref{thm:clt_max_cov_mix_H0}.

As a global generalization of $H_{1}^{L}$, we may write $H_{1}$ in discrete
form as 
\begin{equation}
H_{1}:E[X_{t}X_{t+h}]=\gamma _{h}+c_{h}(t/T),  \label{H1}
\end{equation}%
where as above $c_{h}(\cdot )$ satisfies (\ref{H1_gam(u)}). In this case $%
g_{0}$ $\equiv $ $\lim_{T\rightarrow \infty }T^{-1}\sum_{t=1}^{T}E[X_{t}^{2}]
$ is identically:%
\begin{equation*}
g_{0}=\gamma _{0}+\lim_{T\rightarrow \infty }\frac{1}{T}%
\sum_{t=1}^{T}c_{0}(t/T)=\gamma _{0}+\int_{0}^{1}c_{0}\left( u\right) du>0.
\end{equation*}%
Repeating the above derivations, we find similar to (\ref{EXXW_HL1})%
\begin{equation*}
\frac{1}{T}\sum_{t=1}^{T-h}E\left[ X_{t}X_{t+h}\right] B_{k}(t)=\gamma _{h}%
\frac{1}{T}\sum_{t=1}^{T-h}B_{k}(t)+\frac{1}{T}%
\sum_{t=1}^{T-h}c_{h}(t/T)B_{k}(t)=\int_{0}^{1}c_{h}\left( u\right) \mathcal{%
B}_{k}(u)du+o(1),
\end{equation*}%
and therefore%
\begin{eqnarray*}
\sqrt{T}\left( \hat{\rho}_{h}^{(k)}-\hat{\rho}_{h}\right)  &=&\frac{1}{g_{0}}%
\frac{1}{\sqrt{T}}\sum_{t=1}^{T-h}\left( X_{t}X_{t+h}-E\left[ X_{t}X_{t+h}%
\right] \right) B_{k}(t) \\
&&+\frac{1}{g_{0}}\sqrt{T}\left( \int_{0}^{1}c_{h}\left( u\right) \mathcal{B}%
_{k}(u)du+o(1)\right) +O_{p}(1/\sqrt{T}).
\end{eqnarray*}%
Thus $\max_{\mathcal{H}_{T},\mathcal{K}_{T}}$ $|\sqrt{T}(\hat{\rho}_{h}^{(k)}
$ $-$ $\hat{\rho}_{h})|$ $\overset{p}{\rightarrow }$ $\infty $ given $\lim
\inf_{T\rightarrow \infty }\max_{\mathcal{H}_{T},\mathcal{K}%
_{T}}|\int_{0}^{1}c_{h}\left( u\right) \mathcal{B}_{k}(u)du|$ $>$ $0$.
Identical arguments hold under physical dependence by replacing Lemma \ref%
{lm:clt_max_cov_mix} and Theorem \ref{thm:clt_max_cov_mix_H0} with Theorem %
\ref{thm:clt_max_cov_phys_H0} (cf. Lemma \ref{lm:clt_max_cov_mix}$^{\ast }$
and Theorem \ref{thm:clt_max_cov_mix_H0}$^{\ast }$ in \cite%
{supp_mat_covstat_2024}).

The next result summarizes the preceding discussion.

\begin{theorem}
\label{thm:clt_max_cov_H1}Let Assumption \ref{assum_1abcd}.b,c,d\ hold, and
let $\{\mathcal{H}_{T},\mathcal{K}_{T}\}$ satisfy $0$ $\leq $ $\mathcal{H}%
_{T}$ $\leq $ $T$ $-$ $1$, $\mathcal{H}_{T}$ $\rightarrow $ $\infty $, $%
\mathcal{K}_{T}$ $\rightarrow $ $\infty $, and $\mathcal{H}_{T}$ $=$ $%
O(T^{1/9}(\ln (T))^{1/3})$ under Assumption \ref{assum_1abcd}.a, or $%
\mathcal{H}_{T}$ $=$ $o(T)$ under Assumption \ref{assum_1abcd}.a$^{\ast }$%
.\medskip \newline
$a.$ Under $H_{1}^{L}$, (\ref{max_p_H1}) holds for non-zero $\mathcal{C}%
(h,k) $ in (\ref{C(hk)}), and any sequence $\{\mathcal{K}_{T}\}$ with $%
\mathcal{K}_{T}$ $=$ $o(T^{\kappa })$ for some finite $\kappa $ $>$ $0$ and $%
\eta (\mathcal{K}_{T})$ $=$ $o(\sqrt{T})$.$\medskip $\newline
$b.$ Under $H_{1}$, $\max_{\mathcal{H}_{T},\mathcal{K}_{T}}$ $|\sqrt{T}(\hat{%
\rho}_{h}^{(k)}$ $-$ $\hat{\rho}_{h})|$ $\overset{p}{\rightarrow }$ $\infty $
for any $\{\mathcal{K}_{T}\}$.
\end{theorem}

In the following example we study a simple break in variance in order to
show how the max-test behaves asymptotically.

\begin{example}[Structural Break in Variance]
\label{ex:var}\normalfont Assume covariances do not depend on time:
\linebreak $E[X_{t}X_{t-h}]$ $=$ $\gamma _{h}$ for every $h$ $\geq $ $1$,
and there is a structural break in variance at mid-sample, cf. \cite%
{Perron2006}: 
\begin{equation*}
E[X_{t}^{2}]=g_{1,T}\text{ for }t=1,...,[T/2]\text{ and }E[X_{t}^{2}]=g_{2,T}%
\text{ for }t=[T/2]+1,...,T
\end{equation*}%
for some strictly positive finite sequences $\{g_{1,T},g_{2,T}\}$, $g_{1,T}$ 
$\neq $ $g_{2,T}$. In terms of Walsh or composite Haar systematic samples
and $H_{1}^{L}$, this translates to 
\begin{equation*}
c_{0}(u)=c_{0,1}>0\text{ for }u\in \lbrack 0,1/2)\text{, and }%
c_{0}(u)=c_{0,2}>0\text{ for }u\in \lbrack 1/2,1],
\end{equation*}%
where $c_{0,1}$ $\neq $ $c_{0,2}$. All other $c_{h}\left( u\right) $ $=$ $0$
on $[0,1]$, $h$ $\geq $ $1.$ Hence, by construction of the first Walsh
function $W_{1}(u)$ (or Haar composite $\psi _{1}(u)$): 
\begin{equation*}
\int_{0}^{1}c_{0}\left( u\right) W_{1}(u)du=\int_{0}^{1/2}c_{0}\left(
u\right) du-\int_{1/2}^{1}c_{0}\left( u\right) du=\frac{c_{0,1}-c_{0,2}}{2}%
\neq 0.
\end{equation*}%
Furthermore:\ 
\begin{equation*}
g_{0}\equiv \lim_{T\rightarrow \infty }\frac{1}{T}%
\sum_{t=1}^{T}E[X_{t}^{2}]=\gamma _{0}+\int_{0}^{1}c_{0}\left( u\right)
du=\gamma _{0}+\frac{c_{0,1}+c_{0,2}}{2}.
\end{equation*}%
The normalized correlation difference therefore satisfies for $h$ $\geq $ $1$%
, 
\begin{equation*}
\sqrt{T}\left( \hat{\rho}_{h}^{(k)}-\hat{\rho}_{h}\right) =\frac{1}{g_{0}}%
\frac{1}{\sqrt{T}}\sum_{t=1}^{T-h}\left( X_{t}X_{t+h}-E\left[ X_{t}X_{t+h}%
\right] \right) B_{k}(t)+o_{p}(1).
\end{equation*}%
Under $H_{1}$ at lag $h$ $=$ $0$ and $k$ $=$ $1$ we then have:%
\begin{equation*}
\sqrt{T}\left( \hat{\rho}_{0}^{(1)}-\hat{\rho}_{0}\right) =\frac{1}{g_{0}}%
\frac{1}{\sqrt{T}}\sum_{t=1}^{T}\left( X_{t}^{2}-E\left[ X_{t}^{2}\right]
\right) \mathcal{W}_{1}(t)+\sqrt{T}\left( \frac{c_{0,1}-c_{0,2}}{2g_{0}}%
\right) +o_{p}(1)=\mathcal{Z}_{T}+\mathcal{C}_{T}+o_{p}(1),
\end{equation*}%
say. In view of asymptotic normality of $\mathcal{Z}_{T}$, and $|\mathcal{C}%
_{T}|$ $\rightarrow $ $\infty $, the max-correlation difference test is
consistent when only the variance $E[X_{t}^{2}]$ exhibits a break given $%
\max_{\mathcal{H}_{T},\mathcal{K}_{T}}$ $\sqrt{T}|\hat{\rho}_{h}^{(k)}$ $-$ $%
\hat{\rho}_{h}|$ $\geq $ $\sqrt{T}|\hat{\rho}_{0}^{(1)}$ $-$ $\hat{\rho}_{0}|
$ $=$ $\sqrt{T}|\mathcal{Z}_{T}$ $+$ $\mathcal{C}_{T}|$ $\overset{p}{%
\rightarrow }$ $\infty $.
\end{example}

\section{Dependent wild bootstrap\label{sec:bootstrap}}

We exploit a blockwise wild (multiplier) bootstrap for p-value approximation %
\citep[cf.][]{Liu1988}. The method appears in various places as a multiplier
bootstrap extension of block-based bootstrap methods %
\citep[e.g.][]{Kunsch1989}. \cite{Shao2010} presents a general
nonoverlapping dependent wild bootstrap, exploiting a class of kernel
smoothing weights that omits the truncated kernel, and uses only
\textquotedblleft big\textquotedblright\ blocks of data (\textquotedblleft
little\textquotedblright\ block size is effectively zero). \cite%
{Shao2011_JoE} uses the same method exclusively with a truncated kernel for
a white noise test for a stationary process that is a measurable function of
an iid sequence. In both cases a sequence $\{X_{t}\}_{t=1}^{T}$ is
decomposed into $[T/b_{T}]$ blocks of size $1$ $\leq $ $b_{T}$ $<$ $T$, $%
b_{T}$ $\rightarrow $ $\infty $ and $b_{T}$ $=$ $o(T)$.

\cite{Chernozhukov_etal2019} exploit a Bernstein-like \textquotedblleft
big\textquotedblright\ and \textquotedblleft little\textquotedblright\ block
multiplier bootstrap for high dimensional sample means of stationary,
dependent and bounded sequences. They apply a wild bootstrap on big blocks
and effectively remove the asymptotically negligible little blocks. \cite%
{ZhangWu2014} expand that method for stationary processes by using two
mutually independent iid sequences, one each for big and small blocks.

We expand ideas in \cite{Shao2011_JoE} to non-stationary sequences. The use
of only one set of \textquotedblleft big\textquotedblright\ blocks and a
truncated kernel eases technical arguments and notation, but a more general
use of smoothing kernels and big/little blocks is readily supported by the
theory presented here.

Set a block size $b_{T}$ such that $1\leq b_{T}<T$, $b_{T}/T^{\iota }$ $%
\rightarrow $ $\infty $ and $b_{T}/T^{1-\iota }$ $\rightarrow $ $0$ for some
tiny $\iota $ $>$ $0$. The number of blocks is $\mathcal{N}_{T}$ $=$ $%
[T/b_{T}]$. Denote index blocks $\mathfrak{B}_{s}=\{(s-1)b_{T}+1,\dots
,sb_{T}\}$ for $s=1,\dots ,\mathcal{N}_{T}$, and $\mathfrak{B}_{\mathcal{N}%
_{T}+1}$ $=$ $\{\mathcal{N}_{T}b_{T},...,T\}$. Generate iid random numbers $%
\{\xi _{1},\dots ,\xi _{\mathcal{N}_{T}}\}$ with $E[\xi _{i}]$ $=$ $0$, $%
E[\xi _{i}^{2}]$ $=$ $1$, and $E[\xi _{i}^{4}]$ $<$ $\infty $. Typically $%
\xi _{i}$ is iid $N(0,1)$ in practice, and we make that assumption here to
shorten proofs of weak convergence, cf. Lemmas B.4 and B.4$^{\ast }$ in \cite%
{supp_mat_covstat_2024}. See the proof of Lemma B.4 for further comments on $%
\xi _{i}$.

Define an auxiliary variable $\varphi _{t}=\xi _{s}$ if $t$ $\in $ $%
\mathfrak{B}_{s}$, and let $\Delta \hat{g}_{T}^{(dw)}(h,k)$ be a (centered)
bootstrapped version of $\hat{\gamma}_{h}^{(k)}$ $-$ $\hat{\gamma}_{h}$ $=$ $%
T^{-1}\sum_{t=1}^{T-h}X_{t}X_{t+h}B_{k}(t)$:%
\begin{equation}
\Delta \hat{g}_{T}^{(dw)}(h,k)\equiv \frac{1}{T}\sum_{t=1}^{T-h}\varphi
_{t}\left\{ X_{t}X_{t+h}B_{k}(t)-\frac{1}{T}%
\sum_{s=1}^{T-h}X_{s}X_{s+h}B_{k}(s)\right\} .  \label{DgT}
\end{equation}%
An asymptotically equivalent technique centers only on $X_{t}X_{t+h}$, the
key stochastic term: 
\begin{equation*}
\Delta \hat{g}_{T}^{(dw)}(h,k)\equiv \frac{1}{T}\sum_{t=1}^{T-h}\varphi
_{t}\left\{ X_{t}X_{t+h}-\frac{1}{T}\sum_{s=1}^{T-h}X_{s}X_{s+h}\right\}
B_{k}(t).
\end{equation*}

The bootstrapped test statistic is then $\mathcal{M}_{T}^{(dw)}$ $\equiv $ $%
\hat{\gamma}_{0}^{-1}\max_{\mathcal{H}_{T},\mathcal{K}_{T}}|\sqrt{T}\Delta 
\hat{g}_{T}^{(dw)}(h,k)|.$ Repeat $M$ times. As a by-product of the main
result below, conditional on the sample $\{X_{t}\}_{t=1}^{T}$ this results
in a sequence $\{\mathcal{M}_{T,i}^{(dw)}\}_{i=1}^{M}$ of iid draws $%
\mathcal{M}_{T,i}^{(dw)}$ from the limit null distribution of $\mathcal{M}%
_{T}$ as $T$ $\rightarrow $ $\infty $ asymptotically with probability
approaching one. The approximate p-value is $\hat{p}_{T,M}^{(dw)}$ $\equiv $ 
$1/M\sum_{i=1}^{M}I(\mathcal{M}_{T,i}^{(dw)}$ $\geq $ $\mathcal{M}_{T})$.
The bootstrap test rejects $H_{0}$ at significance level $\alpha $ when $%
\hat{p}_{T,M}^{(dw)}$ $<$ $\alpha $.

The multiplier bootstrap has been studied in many contexts. Consult, e.g., 
\cite{Liu1988}, \cite{Shao2010}, and \cite{Shao2011_JoE} to name a few.
Centering $X_{t}X_{t+h}B_{k}(t)$ $-$ $1/T\sum_{s=1}^{T-h}X_{s}X_{s+h}B_{k}(s)
$ is required because we use $\{X_{t}X_{t+h}B_{k}(t)\}$ to approximate the
null distribution, \textit{whether it is true or not}, and $%
E[X_{t}X_{t+h}B_{k}(t)]$ $\neq $ $0$ for some $(h,k)$ under $H_{1}$. The
block-wise independent zero-mean Gaussian multiplier $\varphi _{t}$ serves
the purpose that $\varphi _{t}\{X_{t}X_{t+h}B_{k}(t)$ $-$ $%
1/T\sum_{s=1}^{T-h}X_{s}X_{s+h}B_{k}(s)\}$, conditioned on the sample $%
\mathfrak{X}_{T}$ $\equiv $ $\{X_{t}\}_{t=1}^{T}$, is zero mean normally
distributed; indeed, $\Delta \hat{g}_{T}^{(dw)}(h,k)|\mathfrak{X}_{T}$ $\sim 
$ $N(0,\mathcal{V}_{T}(h,k))$ for some $\mathcal{V}_{T}(h,k)$ $>$ $0$ $a.s$.
The blocks are constructed such that the dispersion term $\mathcal{V}%
_{T}(h,k)$\ well approximates the null limiting variance under general
dependence, that is $\mathcal{V}_{T}(h,k)$ $\overset{p}{\rightarrow }$ $%
\sigma ^{2}(h,k)$. Thus, in the jargon of \citet[Section
3]{GineZinn1990}, $\Delta \hat{g}_{T}^{(dw)}(h,k)|\mathfrak{X}_{T}$ $\overset%
{d}{\rightarrow }$ $\boldsymbol{Z}(h,k)$ $\sim $ $N(0,\sigma ^{2}(h,k))$ 
\textit{in probability}, ensuring the bootstrapped process yields the null
distribution, irrespective of whether $H_{0}$ holds or not.

Recall $z_{t}(h,k)$ in (\ref{zthk}) and $\mathcal{Z}_{T}(h,k)\equiv 1/\sqrt{T%
}\sum_{t=1}^{T-h}z_{t}(h,k)$. Write $\dddot{g}_{T}(h,k)$ $\equiv $ $1/(T$ $-$
$h)\sum_{u=1}^{T-h}$ \linebreak $E[X_{u}X_{u+h}]B_{k}(u)$ and 
\begin{equation*}
\mathfrak{X}_{T,l}(h,k)\equiv \sum_{t=(l-1)b_{T}+1}^{lb_{T}}\left\{
X_{t}X_{t+h}B_{k}(t)-\dddot{g}_{T}(h,k)\right\} ,
\end{equation*}%
and define pre-asymptotic and asymptotic long run covariance functions $%
s_{T}^{2}(h,k;\tilde{h},\tilde{k})$ $\equiv $ \linebreak $%
1/T\sum_{l=1}^{(T-h\vee \tilde{h})/b_{T}}E[\mathfrak{X}_{T,l}(h,k)\mathfrak{X%
}_{T,l}(\tilde{h},\tilde{k}))]$ and $s^{2}(h,k;\tilde{h},\tilde{k})$ $\equiv 
$ $\lim_{T\rightarrow \infty }s_{T}^{2}(h,k;\tilde{h},\tilde{k})$.

\begin{assumption}
\label{assum_boot} \ \ \medskip \newline
$a.$ $(i)$ $\lim \inf_{T\rightarrow \infty }s_{T}^{2}(h,k;\tilde{h},\tilde{k}%
)$ $>$ $0$ $\forall (h,\tilde{h},k,\tilde{k})$; $(ii)$ $\max_{\mathcal{H}%
_{T},\mathcal{K}_{T}}|s_{T}^{2}(h,k;\tilde{h},\tilde{k})$ $-$ $s^{2}(h,k;%
\tilde{h},\tilde{k})|$ $=$ $O(T^{-\iota })$ for some infinitessimal $\iota $ 
$>$ $0$.\medskip \newline
$b.$ $b_{T}/T^{\iota }$ $\rightarrow $ $\infty $ and $b_{T}$ $=$ $%
o(T^{1/2-\iota })$ for some infinitessimal $\iota $ $>$ $0$.
\end{assumption}

\begin{remark}
\normalfont($a.i$) is the fourth order block bootstrap version of Assumption %
\ref{assum_1abcd}.c, used to ensure a high dimensional central limit theory
extends to a long run bootstrap variance, cf. 
\citet[Lemma
3.1]{Chernozhukov_etal2013}. ($a.ii$) seems unavoidable, and is required to
link covariance functions for a high dimensional bootstrap theory, cf. %
\citet[Lemma 3.1]{Chernozhukov_etal2013} and 
\citet[Theorem 2, Proposition
1]{Chernozhukov_etal2015}. The property is trivial under stationary
geometric mixing or physical dependence, and otherwise restricts the degree
of allowed heterogeneity. ($b$) simplifies a bootstrap weak convergence
proof, but can be weakened at the cost of added notation, e.g. $b_{T}/(\ln
(T))^{a}$ $\rightarrow $ $\infty $ and $b_{T}$ $=$ $o(T^{1/2}/(\ln (T))^{b})$
for some $a,b$ $>$ $0$.
\end{remark}

The blockwise wild bootstrap is valid asymptotically under mixing or
physical dependence.

\begin{theorem}
\label{thm:p_dep_wild_boot_mix}Let Assumptions \ref{assum_1abcd}.b,c,d and %
\ref{assum_boot} hold, let $\mathcal{H}_{T},\mathcal{K}_{T}$ $\rightarrow $ $%
\infty $, and let the number of bootstrap samples $M$ $=$ $M_{T}$ $%
\rightarrow $ $\infty $ as $T$ $\rightarrow $ $\infty $. Let $\{b_{T},%
\mathcal{H}_{T}\}$ satisfy $b_{T}$ $\rightarrow $ $\infty $ and $b_{T}$ $=$ $%
O(T^{1/2-\iota })$, $0$ $\leq $ $\mathcal{H}_{T}$ $\leq $ $T$ $-$ $1$, and
under Assumption \ref{assum_1abcd}.a $\mathcal{H}_{T}$ $=$ $O(T^{1/9}(\ln
(T))^{1/3})$, or $\mathcal{H}_{T}$ $=$ $O(T^{1-\iota }/b_{T})$ for tiny $%
\iota $ $>$ $0$\ under Assumption \ref{assum_1abcd}.a$^{\ast }$. Under $%
H_{0} $, $P(\hat{p}_{T,M}^{(dw)}$ $<$ $\alpha )$ $\rightarrow $ $\alpha $
for any sequence $\{\mathcal{K}_{T}\}$ satisfying $\mathcal{K}_{T}$ $=$ $%
o(T^{\kappa })$ for some finite $\kappa $ $>$ $0$ and $\eta (\mathcal{K}%
_{T}) $ $=$ $o(\sqrt{T})$. Under $H_{1}$ in (\ref{H1}) where $c_{h}(\cdot )$
satisfy (\ref{cW}), $P(\hat{p}_{T,M}^{(dw)}$ $<$ $\alpha )$ $\rightarrow $ $%
1 $ for any $\{\mathcal{K}_{T}\}$.
\end{theorem}

\section{Monte carlo study\label{sec:monte_carlo}}

We now study the proposed bootstrap test in a controlled environment. We
generate $1000$ independently drawn samples from various models, with sample
sizes $T$ $\in $ $\{64,128,256,512\}$. The models under the null and
alternative hypotheses are detailed below.

\subsection{Empirical size}

We use four models of covariance stationary processes: MA(1), AR(1), Self
Exciting Threshold AR(1) [SETAR], and GARCH (1,1):

\begin{center}
\begin{tabular}{lll}
null-1 & MA(1) & $X_{t}=\epsilon _{t}$ \\ 
null-2 & AR(1) & $X_{t}=.5X_{t-1}+\epsilon _{t}$ \\ 
null-3 & SETAR & $X_{t}=.7X_{t-1}-1.4X_{t-1}I\left( X_{t-1}>0\right)
+\epsilon _{t}\text{,}$ \\ 
null-4 & GARCH(1,1) & $X_{t}=\sigma _{t}z_{t}\text{, }z_{t}\overset{iid}{%
\sim }N(0,1),\text{ }\sigma _{t}^{2}=1+.3\epsilon _{t-1}^{2}+.6\sigma
_{t-1}^{2}$%
\end{tabular}
\end{center}

Models \#1-\#3 have an iid error $\epsilon _{t}$ distributed $N(0,1)$ or
Student's-$t$ with $5$ degrees of freedom ($t_{5}$); or $\epsilon _{t}$ is
stationary GARCH(1,1) $\epsilon _{t}$ $=$ $\sigma _{t}z_{t},$ $z_{t}$ $%
\overset{iid}{\sim }$ $N(0,1)$, $\sigma _{t}^{2}$ $=$ $1$ $+$ $.3\epsilon
_{t-1}^{2}$ $+$ $.6\sigma _{t-1}^{2}$, with iteration $\sigma _{1}^{2}$ $=$ $%
1$ and $\sigma _{t}^{2}$ $=$ $1$ $+$ $.3\epsilon _{t-1}^{2}$ $+$ $.6\sigma
_{t-1}^{2}$ for $t$ $=$ $2,...,T$. \ The SETAR model switches between AR(1)
regimes with correlations $.7$ and $-.7$. GARCH and SETAR models, and any
model with GARCH errors, do not have a linear form $X_{t}$ $=$ $%
\sum_{i=0}^{\infty }\psi _{i}Z_{t-i}$, with iid $Z_{t}$ and non-random $\psi
_{i}$. We simulate $2T$ observations for each model and retain the latter $T$
observations for analysis. Test results in GARCH cases should be viewed with
caution: max-test asymptotics have only been established under
sub-exponentail tails while GARCH processes have regularly varying (i.e.
\textquotedblleft power\textquotedblright ) tails, and JWW's test requires a
linear model with an iid error.

\subsection{Empirical power}

We study empirical power by using models similar to those in \cite%
{Paparoditis_2010}; \cite{Dette_etal_2011}, \cite{PreussVetterDette2013} and
JWW, with the addition of allowing for non-iid errors and non-stationarity
in variance. The models are as follows:

\begin{center}
\begin{tabular}{lll}
alt-1 & (NI) & $X_{t}=1.1\cos \{1.5-\cos (4\pi t/T)\}\epsilon
_{t-1}+\epsilon _{t}$ \\ 
alt-2 & (NVIII) & $X_{t}=.8\cos \{1.5-\cos (4\pi t/T)\}\epsilon
_{t-6}+\epsilon _{t}$ \\ 
alt-3 & (NII) & $X_{t}=.6\times \sin (4\pi t/T)X_{t-1}+\epsilon _{t}$ \\ 
alt-4 & $\text{(NIII)}$ & $X_{t}=\left\{ 
\begin{array}{cl}
.5X_{t-1}+\epsilon _{t} & \text{for }\{1\leq t\leq T/4\}\cup \{3T/4<t\leq T\}
\\ 
-.5X_{t-1}+\epsilon _{t} & \text{for }T/4<t\leq 3T/4%
\end{array}%
\right. $%
\end{tabular}

\begin{tabular}{lll}
alt-5 & $\text{(NVI)}$ & $X_{t}=\left\{ 
\begin{array}{cl}
.5X_{t-1}+\epsilon _{t} & \text{for }1\leq t\leq T/2 \\ 
-.5X_{t-1}+\epsilon _{t} & \text{for }T/2<t\leq T%
\end{array}%
\right. $ \\ 
alt-6 & $\text{(eq. (16))}$ & $X_{t}=2\epsilon _{t}-\left\{ 1+.5\cos (2\pi
t/T)\right\} \epsilon _{t-1}$ \\ 
alt-7 & $\text{(NV)}$ & $X_{t}=-.9\sqrt{(t/T)}X_{t-1}+\epsilon _{t}$ \\ 
alt-8 &  & $X_{t}=.5X_{t-1}+v_{t}\text{: }\left\{ 
\begin{array}{ll}
v_{t}=\epsilon _{t} & \text{for }1\leq t\leq 3T/4 \\ 
v_{t}=2\epsilon _{t} & \text{for }3T/4<t\leq T%
\end{array}%
\right. $ \ \ \ \ \ \ \ \ \ \ \ \ \ \ \  \\ 
alt-9 &  & $X_{t}=.8\cos \{1.5-\cos (4\pi t/T)\}\epsilon _{t-25}+\epsilon
_{t}$%
\end{tabular}
\end{center}

Models 1-7 are used in JWW: we display parenthetically their corresponding
model/equation number. Models 1, 2, 4 are considered in \cite%
{Paparoditis_2010}; \cite{Dette_etal_2011} use models 1, 2, 4, and 6; \cite%
{PreussVetterDette2013} study 2, 5, and 7. Alt-8 presents a structural
change in variance only, and alt-9 is a distant version of alt-2 and
therefore more difficult to detect (lag $25$ as opposed to lag $6$). As
above, we use either iid standard normal, iid $t_{5}$, or GARCH(1,1) $%
\epsilon _{t}$.

\subsection{Tests}

\subsubsection{Max-test}

We perform the bootstrapped max-correlation difference test with $\mathcal{M}%
_{T}$ and $\mathcal{M}_{T}^{(p)}$. The latter has penalties $p_{h}$ $=$ $(h$ 
$+$ $1)^{1/4}/2$ and $q_{k}$ $=$ $k^{1/4}/2$. More severe penalties, e.g. $%
q_{k}$ $=$ $k^{1/2}/2$, do not improve test performance. A weighted version
of the test with HAC estimator (\ref{V_hac}) leads to competitive size but
generally lower power, hence we focus only on $\mathcal{M}_{T}$ and $%
\mathcal{M}_{T}^{(p)}$. We use Walsh or Haar functions for two max-tests,
and a third combined max-max-statistic shown below (\ref{M_combo}). We only
report results based on Walsh functions because ($i$) the Haar-based tests
(max-test, and JWW's test detailed below) yielded far lower power across
most alternatives studied here; hence ($ii$) the max-max test performed
essentially on par with, or was slightly trumped by, the Walsh-based test.

We use $500$ bootstrap samples with multiplier iid variable $\xi _{t}$ $\sim 
$ $N(0,1)$. Theorem \ref{thm:p_dep_wild_boot_mix} requires a block size
bound $b_{T}$ $=$ $o(T^{1/2-\iota })$ for some tiny $\iota $ $>$ $0$, hence
we use $b_{T}$ $=$ $[T^{1/2-\eta }]$ where $\eta $ $=$ $10^{-10}$. Similar
block sizes, e.g. $b_{T}$ $=$ $[bT^{1/2-\eta }]$ with $b$ $\in \lbrack .5,2]$
lead to similar results.\footnote{\cite{Shao2011_JoE} uses $b_{T}$ $=$ $%
[bT^{1/2}]$ with $b$ $\in $ $\{.5,1,2\}$, leading to qualitatively similar
results. \cite{HillMotegi2020} also use $b$ $=$ $1$, but find qualitatively
similar results for values $b$ $\in $ $\{.5,1,2\}$.}

Theorem \ref{thm:p_dep_wild_boot_mix} also requires $\mathcal{H}_{T}$ $=$ $%
O(T^{1-\iota }/b_{T})$, $\mathcal{K}_{T}$ $=$ $o(T^{\kappa })$ for some $%
\kappa $ $>$ $0$, and $\eta (\mathcal{K}_{T})$ $=$ $o(\sqrt{T})$. In the
Walsh case $\eta (k)$ $=$ $k$ hence $\mathcal{K}_{T}$ $=$ $o(\sqrt{T})$; in
the Haar case $\eta (k)$ $=$ $2^{k}$ hence $\mathcal{K}_{T}$ $=$ $o(\ln (T))$%
. In the Walsh case, we used two pairings of sequences $\{\mathcal{H}_{T},%
\mathcal{K}_{T}\}$. The first $\mathcal{H}_{T}$ $=$ $[\log _{2}(T)^{.99}$ $-$
$3]$ and $\mathcal{K}_{T}$ $=$ $[T^{1/3}]$ is used in JWW. The second $%
\mathcal{H}_{T}$ $=$ $[2T^{.49}]$ and $\mathcal{K}_{T}$ $=$ $[.5T^{.49}]$
satisfies our assumptions but are not valid in JWW. The latter $(\mathcal{H}%
_{T},\mathcal{K}_{T})$ are generally larger, where $\mathcal{H}_{T}$ is
larger by an order of $\times 7$. This will lead to higher power for large $%
T $ in theory, but in small samples obviously\ a larger $h$ results in fewer
observations for computation, and therefore a loss in sharpness in
probability. In the Haar case we use either $\mathcal{H}_{T}$ above, and $%
\mathcal{K}_{T}$ $=$ $[(\ln (T))^{.99}]$. Refer to Table \ref{tbl:HtKt}.

\begin{table}[tbp]
\caption{$\mathcal{H}_{T}$, $\mathcal{K}_{T}$ Combinations by Basis}
\label{tbl:HtKt}\centering%
\begin{tabular}{l|cc|cc|cc|c|c|}
& \multicolumn{4}{|c|}{Walsh Basis $\{W_{k}(t)\}$} & \multicolumn{4}{|c|}{
Haar Basis $\{\Psi _{k}(t)\}$} \\ \hline\hline
& \multicolumn{2}{|c|}{Case 1 (JWW)} & \multicolumn{2}{|c|}{Case 2} & 
\multicolumn{2}{|c}{Case 1 (JWW)} & \multicolumn{2}{|c|}{Case 2} \\ 
\hline\hline
$T$ & $H_{T}$ & $K_{T}$ & $H_{T}$ & $K_{T}$ & $H_{T}$ & $K_{T}$ & $H_{T}$ & $%
K_{T}$ \\ \hline\hline
& $\log _{2}(T)^{.99}$ $-$ $3$ & $T^{1/3}$ & $2T^{.49}$ & $.5T^{.49}$ & $%
\log _{2}(T)^{.99}$ $-$ $3$ & \multicolumn{1}{c|}{$(\ln (T))^{.99}$} & $%
2T^{.49}$ & $(\ln (T))^{.99}$ \\ \hline
$64$ & $2$ & $4$ & $14$ & $3$ & $2$ & \multicolumn{1}{c|}{$4$} & $14$ & $4$
\\ 
$128$ & $3$ & $5$ & $20$ & $5$ & $3$ & \multicolumn{1}{c|}{$5$} & $20$ & $5$
\\ 
$256$ & $4$ & $6$ & $30$ & $7$ & $4$ & \multicolumn{1}{c|}{$5$} & $30$ & $5$
\\ 
$512$ & $5$ & $8$ & $42$ & $10$ & $5$ & \multicolumn{1}{c|}{$6$} & $42$ & $6$
\\ \hline\hline
\end{tabular}%
\end{table}

\subsubsection{JWW test}

Write $\boldsymbol{\hat{\gamma}}_{h}$ $\equiv $ $[\hat{\gamma}_{1},...,\hat{%
\gamma}_{h}]^{\prime }$ and $\boldsymbol{\hat{\gamma}}_{h}^{(k)}$ $\equiv $ $%
[\hat{\gamma}_{1}^{(k)},...,\hat{\gamma}_{h}^{(k)}]^{\prime }$. The test
statistic is%
\begin{equation*}
\mathcal{\hat{D}}_{T}\equiv \max_{1\leq k\leq \mathcal{K}_{T}}\left[
\max_{1\leq h\leq \mathcal{H}_{T}}\left\{ T\left( \boldsymbol{\hat{\gamma}}%
_{h}^{(k)}-\boldsymbol{\hat{\gamma}}_{h}\right) ^{\prime }\left( \hat{\Gamma}%
_{h}^{(k)}\right) ^{-1}\left( \boldsymbol{\hat{\gamma}}_{h}^{(k)}-%
\boldsymbol{\hat{\gamma}}_{h}\right) -2h\right\} -\sqrt{k-1}\right] 
\end{equation*}%
where $\hat{\Gamma}_{h}^{(k)}$ is an estimator of the $h\times h$ asymptotic
covariance matrix of $\sqrt{T}(\boldsymbol{\hat{\gamma}}_{h}^{(k)}-%
\boldsymbol{\hat{\gamma}}_{h})$. See 
\citet[Sections
2.3-2.5]{JinWangWang2015} for details on computing $\hat{\Gamma}_{h}^{(k)}$
(under the assumption of linearity $X_{t}$ $=$ $\sum_{i=0}^{\infty }\psi
_{i}Z_{t-i}$ with an iid $Z_{t}$).\footnote{%
There is a typo in \citet[Theorem 2]{JinWangWang2015} concerning their
covariance matrix and therefore its estimator. A parameter $\kappa _{4}$,
referred to as the kurtosis of the iid $Z_{t}$, is in fact the excess
kurtosis (\textit{kurtosis} $-3$). See Proposition 7.3.1 in \cite%
{BrockwellDavis2991}, in particular eq. (7.3.5), cf. 
\citet[p.
915]{JinWangWang2015}.} We use both Walsh and Haar bases, the same tuning
parameters that JWW use for covariance matrix estimation, and the same $\{%
\mathcal{H}_{T},\mathcal{K}_{T}\}$ described above.\footnote{%
The bandwidth parameter $\lambda $ in $[T^{\lambda }]$, the number of sample
covariances that enter the asymptote covariance matrix estimator, is set to $%
\lambda $ $=$ $.4$ based on a private communication with the authors. In
order to compute the (excess) kurtosis of iid $Z_{t}$ under linearity,
similar to \citet[eq. (15)]{JinWangWang2015} we use an estimator in \cite%
{KreissPaparoditis2015}, with two bandwidths $b_{j}$ $=$ $c_{j}T^{-1/3}$
where each $c_{j}$ $=1.25\times \hat{\gamma}(0)$ \citep[see][p.
903]{JinWangWang2015}.}

We perform the test both based on a simulated critical values (denoted $%
\mathcal{\hat{D}}_{T}^{cv}$), and bootstrapped p-values ($\mathcal{\hat{D}}%
_{T}^{dw}$) in order to make a direct comparison with the method developed
here. We simulate critical values for each basis and each pair $(\mathcal{H}%
_{T},\mathcal{K}_{T})$ by running a separate simulation with $200,000$
independently drawn samples of size $T$ of iid $N(0,1)$ distributed random
variables $X_{t}$, and use the true excess kurtosis value $0$ in the
covariance estimator $\hat{\Gamma}_{h}^{(k)}$. The bootstrap is performed by
replacing $\boldsymbol{\hat{\gamma}}_{h}^{(k)}-\boldsymbol{\hat{\gamma}}_{h}$
in $\mathcal{\hat{D}}_{T}$ with $\Delta \hat{g}_{T}^{(dw)}(h,k)$ from (\ref%
{DgT}). We do not prove asymptotic validity of the bootstrapped p-value, but
once uniform consistency of $\hat{\Gamma}_{h}^{(k)}$ is established, it
follows identically from arguments given in the proof of Theorem \ref%
{thm:p_dep_wild_boot_mix}. Indeed, the bootstrap is valid for linear and
nonlinear processes with iid or non-iid innovations, and covering the
nonstationary processes under $H_{1}$. The simulated critical values,
however, are suitable in theory only for linear processes with iid
innovations since they rely on the specific form of $\hat{\Gamma}_{h}^{(k)}$
used here, and a pivotal Gaussian null limit distribution, cf. 
\citet[Sections
2.3-2.5]{JinWangWang2015}.

\subsection{Results}

Tables A.3-A.6 in \citet[Appendix D]{supp_mat_covstat_2024} present
rejection frequencies at\ $(1\%,5\%,10\%)$ significance levels when a Walsh
basis is used.

The penalized max-test does not perform better than the non-penalized test,
and generally performs worse under the alternative. Indeed, as discussed
above, there is no theory driven reason for adding penalties for a max-test.
In the sequel we therefore only discuss the non-penalized test.

Similarly, the bootstrapped JWW test is generally over-sized, and massively
over-sized at small $T$ under $(\mathcal{H},\mathcal{K})$ Case 1, the only
valid case in this study. We suspect the cause is the estimated variance
matrix due to its many components and tuning parameters. We henceforth only
discuss results based on simulated critical values.

\subsubsection{Null}

Both tests are comparable for MA and AR models with iid Gaussian or $t_{5}$
errors, with fairly accurate empirical size. The max-test has accurate size
in many cases, and is otherwise conservative. JWW's test tends to be
over-sized in the AR model with GARCH errors under both $(\mathcal{H},%
\mathcal{K})$ cases, and is over-sized in the AR model with $t_{5}$ errors
under Case 2 when $T$ $\leq $ $128$. Recall $\mathcal{H}_{T}$ is much larger
under Case 2, which will be a hindrance at smaller $T$ for test statistics
that simultaneously incorporate a set of autocovariances (e.g. Wald or
portmanteau statistics).

In the SETAR case JWW's test is largely over-sized, while the max-test is
slightly under-sized with improvement under $(\mathcal{H},\mathcal{K})$ Case
2. JWW's test is over-sized for small $T$ with the GARCH model, but
otherwise works well.

\subsubsection{Alternative}

In Table \ref{tbl:dominance} we give a simple summary of which test
generally dominates for each model and case based on the complete simulation
results. In brief, each test dominates for certain models, and in some cases
they are comparable. JWW's test generally dominates in models 1, 3, and 4,
and for model 7 for larger sample sizes. This applies across error cases,
including GARCH errors.

The max-test dominates in models 2, 6, 8 and 9, with strong domination for
model 8 (break in variance), and models 2 and 9 (distant nonstationarity).
Indeed, JWW's test has only negligible power for models 2, 8 and 9: by
construction it cannot detect a break in variance (model 8), and seems
incapable of detecting a distant (model 9), or even semi-distant (model 2),
form of covariance nonstationarity.

Overall, both tests clearly have merit, and seem to complement each other
based on the different cases in which they each excel. See 
\citet[Appendix
C]{supp_mat_covstat_2024} for an application of both tests to international
exchange rates.

\begin{table}[tbp]
\caption{Test Dominance Summary}
\label{tbl:dominance}
\begin{center}
\begin{tabular}{l|ccc|ccc}
& \multicolumn{3}{c|}{$(\mathcal{H},\mathcal{K)}$ Case 1} & 
\multicolumn{3}{|c}{$(\mathcal{H},\mathcal{K)}$ Case 2} \\ \hline\hline
$H_{1}\backslash \epsilon _{t}$ & $N(0,1)$ & $t_{5}$ & GARCH & $N(0,1)$ & $%
t_{5}$ & GARCH \\ \hline
\multicolumn{1}{c|}{alt-1} & \multicolumn{1}{|l}{$\mathcal{\hat{D}}_{T}$
small $n$} & \multicolumn{1}{l}{$\mathcal{\hat{D}}_{T}$ small $n$} & 
\multicolumn{1}{l|}{$\mathcal{\hat{D}}_{T}$} & \multicolumn{1}{|l}{$\mathcal{%
\hat{D}}_{T}$ small $n$} & \multicolumn{1}{l}{$\mathcal{\hat{D}}_{T}$ small $%
n$} & \multicolumn{1}{l}{$\mathcal{\hat{D}}_{T}$} \\ 
\multicolumn{1}{c|}{alt-2} & \multicolumn{1}{|l}{\textcolor{red}{%
\boldsymbol{$\mathcal{\hat{M}}_{T}$}}} & \multicolumn{1}{l}{%
\textcolor{red}{\boldsymbol{$\mathcal{\hat{M}}_{T}$}}} & \multicolumn{1}{l|}{%
\textcolor{red}{\boldsymbol{$\mathcal{\hat{M}}_{T}$}}} & \multicolumn{1}{|l}{%
\textcolor{red}{\boldsymbol{$\mathcal{\hat{M}}_{T}$}}} & \multicolumn{1}{l}{%
\textcolor{red}{\boldsymbol{$\mathcal{\hat{M}}_{T}$}}} & \multicolumn{1}{l}{%
\textcolor{red}{\boldsymbol{$\mathcal{\hat{M}}_{T}$}}} \\ 
\multicolumn{1}{c|}{alt-3} & \multicolumn{1}{|l}{$\mathcal{\hat{D}}_{T}$
small $n$} & \multicolumn{1}{l}{$\mathcal{\hat{D}}_{T}$ small $n$} & 
\multicolumn{1}{l|}{$\mathcal{\hat{D}}_{T}$} & \multicolumn{1}{|l}{$\mathcal{%
\hat{D}}_{T}$ small $n$} & \multicolumn{1}{l}{$\mathcal{\hat{D}}_{T}$ small $%
n$} & \multicolumn{1}{l}{$\mathcal{\hat{D}}_{T}$} \\ 
\multicolumn{1}{c|}{alt-4} & \multicolumn{1}{|l}{$\mathcal{\hat{D}}_{T}$
small $n$} & \multicolumn{1}{l}{$\mathcal{\hat{D}}_{T}$ small $n$} & 
\multicolumn{1}{l|}{$\mathcal{\hat{D}}_{T}$} & \multicolumn{1}{|l}{$\mathcal{%
\hat{D}}_{T}$ small $n$} & \multicolumn{1}{l}{$\mathcal{\hat{D}}_{T}$ small $%
n$} & \multicolumn{1}{l}{$\mathcal{\hat{D}}_{T}$} \\ 
\multicolumn{1}{c|}{alt-5} & \multicolumn{1}{|l}{\textcolor{red}{%
\boldsymbol{$\mathcal{\hat{M}}_{T}$}}} & \multicolumn{1}{l}{%
\textcolor{red}{\boldsymbol{$\mathcal{\hat{M}}_{T}$}}} & \multicolumn{1}{l|}{%
\textcolor{red}{\boldsymbol{$\mathcal{\hat{M}}_{T}$}}} & \multicolumn{1}{|l}{%
\textcolor{red}{\boldsymbol{$\mathcal{\hat{M}}_{T}$}}} & \multicolumn{1}{l}{%
\textcolor{red}{\boldsymbol{$\mathcal{\hat{M}}_{T}$}}} & \multicolumn{1}{l}{%
\textcolor{red}{\boldsymbol{$\mathcal{\hat{M}}_{T}$}}} \\ 
\multicolumn{1}{c|}{alt-6} & \multicolumn{1}{|l}{\textcolor{red}{%
\boldsymbol{$\mathcal{\hat{M}}_{T}$}}} & \multicolumn{1}{l}{%
\textcolor{red}{\boldsymbol{$\mathcal{\hat{M}}_{T}$}}} & \multicolumn{1}{l|}{
similar} & \multicolumn{1}{|l}{\textcolor{red}{\boldsymbol{$\mathcal{%
\hat{M}}_{T}$}}} & \multicolumn{1}{l}{\textcolor{red}{\boldsymbol{$\mathcal{%
\hat{M}}_{T}$}}} & \multicolumn{1}{l}{similar} \\ 
\multicolumn{1}{c|}{alt-7} & \multicolumn{1}{|l}{$\mathcal{\hat{D}}_{T}$
larger $n$} & \multicolumn{1}{l}{$\mathcal{\hat{D}}_{T}$ large n} & 
\multicolumn{1}{l|}{$\mathcal{\hat{D}}_{T}$ large $n$} & \multicolumn{1}{|l}{%
$\mathcal{\hat{D}}_{T}$ large $n$} & \multicolumn{1}{l}{$\mathcal{\hat{D}}%
_{T}$ large $n$} & \multicolumn{1}{l}{$\mathcal{\hat{D}}_{T}$ large $n$} \\ 
\multicolumn{1}{c|}{alt-8} & \multicolumn{1}{|l}{\textcolor{red}{%
\boldsymbol{$\mathcal{\hat{M}}_{T}$}}} & \multicolumn{1}{l}{%
\textcolor{red}{\boldsymbol{$\mathcal{\hat{M}}_{T}$}}} & \multicolumn{1}{l|}{%
\textcolor{red}{\boldsymbol{$\mathcal{\hat{M}}_{T}$}}} & \multicolumn{1}{|l}{%
\textcolor{red}{\boldsymbol{$\mathcal{\hat{M}}_{T}$}}} & \multicolumn{1}{l}{%
\textcolor{red}{\boldsymbol{$\mathcal{\hat{M}}_{T}$}}} & \multicolumn{1}{l}{%
\textcolor{red}{\boldsymbol{$\mathcal{\hat{M}}_{T}$}}} \\ 
\multicolumn{1}{c|}{alt-9} & \multicolumn{1}{|l}{\textcolor{red}{%
\boldsymbol{$\mathcal{\hat{M}}_{T}$}}} & \multicolumn{1}{l}{%
\textcolor{red}{\boldsymbol{$\mathcal{\hat{M}}_{T}$}}} & \multicolumn{1}{l|}{%
\textcolor{red}{\boldsymbol{$\mathcal{\hat{M}}_{T}$}}} & \multicolumn{1}{|l}{%
\textcolor{red}{\boldsymbol{$\mathcal{\hat{M}}_{T}$}}} & \multicolumn{1}{l}{%
\textcolor{red}{\boldsymbol{$\mathcal{\hat{M}}_{T}$}}} & \multicolumn{1}{l}{%
\textcolor{red}{\boldsymbol{$\mathcal{\hat{M}}_{T}$}}} \\ \hline\hline
\end{tabular}%
\end{center}
\par
{\small Each cell dictates which test performed best (in certain cases). For
example ``$\mathcal{\hat{D}}_{T}$ small $n$" implies $\mathcal{\hat{D}}_{T}$
dominates for smaller sample sizes, and for other $n$ the two tests are
comparable. ``$\mathcal{\hat{M}}_{T}$" implies $\mathcal{\hat{M}}_{T}$
dominates across sample sizes.}
\end{table}

\section{Conclusion\label{sec:conclusion}}

We present a max-correlation difference test for testing covariance
stationarity in a general setting that allows for nonlinearity and random
volatility, and heterogeneity under either hypothesis. Our test exploits a
generic orthonormal basis under mild conditions, with Walsh and Haar wavelet
function examples. We do not require estimation of an asymptotic covariance
matrix, our test can detect a break in variance, and we deliver an
asymptotically valid dependent wild bootstrapped p-value. Orthonormal basis
based tests direct power toward alternatives implied by basis-specific
systematic samples. Thus, by combining bases a power improvement may be
achievable. In controlled experiments, however, we find the Walsh basis
yields superior results compared to a composite Haar basis. We leave for
future endeavors the question of whether other bases may yet perform better
than the Walsh basis for a test of covariance stationarity.

Furthermore, the max-test dominates JWWs in some case, while JWW's dominates
in others. The max-test is best capable of delivering sharp empirical size
for a nonlinear process and when errors are non-iid, and is particularly
suited for detecting distant (large lag) forms of covariance
non-stationarity, and a break in variance. The former corroborates findings
in \cite{HillMotegi2020}, who find a max-correlation white noise test
strongly dominates Wald and portmanteau tests when there is a distant
non-zero correlation. We conjecture this will carry over to other
nonstationary models with distant breaks in covariance, but leave this idea
for future consideration.

\setcounter{equation}{0} \renewcommand{\theequation}{{\thesection}.%
\arabic{equation}} \appendix

\begin{acks}[Acknowledgments]
 The authors thank two expert referees, and editor Davy Paindaveine, for comments that led to major improvements of the manuscript.
\end{acks}%
\begin{acks}[Supplementary material]
Additional results, omitted proofs and complete simulation results.
\end{acks}

\bibliographystyle{imsart-nameyear}
\bibliography{refs_cov_stat_test}

\end{document}